\documentclass[11pt]{article}

\usepackage{fancyhdr, amsmath, amssymb, amsfonts, url, subfigure, dsfont, mathrsfs, graphicx, epsfig, subfigure, amsthm, tikz}
    \usepackage{epsfig}
    \usepackage{graphicx}
    \usepackage{color}

\newlength{\guillotine}
\newtheorem{thm}{Theorem}[section]
\newtheorem{cor}[thm]{Corollary}
\newtheorem{lemma}[thm]{Lemma}

\newtheorem{question}[thm]{Objective}

\newtheorem{definition}[thm]{Definition}
\newtheorem{example}[thm]{Example}

\theoremstyle{remark}
\newtheorem{rem}[thm]{Remark}

\AtEndDocument{\bigskip{\footnotesize
  \noindent
  \textsc{O. Jenkinson, School of Mathematical Sciences, Queen Mary University of London, Mile End Road, London, E1 4NS, UK}
  
  \noindent
  \textit{E-mail address}: \texttt{o.jenkinson@qmul.ac.uk}
  \par
  \addvspace{\medskipamount}
  \noindent
  \textsc{M. Pollicott, Department of Mathematics, Warwick University, Coventry, CV4 7AL, UK.}
  
  \noindent
  \textit{E-mail address}: \texttt{masdbl@warwick.ac.uk}
}}

\begin{document}
\date{}

\title{A  dynamical approach to validated numerics}

\author{O. Jenkinson and M. Pollicott}

\date{In memorium, A.B. Katok
\thanks{The second author collaborated with A. Katok, G. Knieper and H. Weiss over 30 years ago on the regularity of topological entropy for Anosov flows.  One of the methods used there was based on the characterisation of the entropy as a pole for the dynamical zeta function, which is a complex function closely related to the determinant central to ideas we pursue here.}
}

\maketitle
\begin{abstract}
We describe a method, using periodic points and determinants, 
 for giving alternative expressions for dynamical quantities (including Lyapunov exponents and Hausdorff dimension of invariant sets) associated to analytic hyperbolic systems.  This leads to validated numerical estimates on their values, \end{abstract}
\section{Introduction}

\subsection{An overview}
There are many well known and important numerical invariants in the context of dynamical systems,
including for example entropy, Lyapunov exponents, and the Hausdorff dimension of invariant sets.
Even in the context of uniformly hyperbolic systems, many of these invariants do not have simple explicit expressions, nor are they easy to estimate. This leads  naturally to the search for useful alternative expressions for each of these quantities which, in particular, 
 lend themselves  to accurate  numerical evaluation.

Our aim in this account is to describe an approach developed over some years,  based on periodic points for (hyperbolic) dynamical systems. 
In particular, in addition to presenting a general overview we  will also
give  a general recipe for converting this information  into formulae and useful numerical estimates  for  a variety of characteristic values, such as those mentioned above.  Perhaps the  most important aspect of this side of the work  is that we obtain rigorous bounds on the errors.  
  Given the data available, the idea  is to minimise the error estimate, subject to the practical constraints imposed by the computational power available.  There is  always some flexibility in the choices, which we can exploit 
  in order to obtain the best error estimate, i.e., the smallest bound on the error.

Typically, the quantities that we can expect to study in this approach are those that can be expressed in terms of the thermodynamic pressure.  The main ingredient in our approach is a complex function 
(called the determinant) which packages together the data on periodic orbits and from which can be 
derived estimates on the pressure, and consequently the quantities of interest.

 In our presentation, no specialist knowledge is required to implement the general result.  However, for completeness we give a fairly complete outline of the proof (with some details on operator theory postponed  to   Appendix A).

\subsection{A historical perspective}
The starting point for our 
approach is the important 
work of
Grothendieck in the 1950s on nuclear operators \cite{grothendieckthesis, grothendieckbullsmf}.
  This extended the classical  theory of trace class operators and Fredholm determinants \cite{fredholm}.  
Although the impact of this theory in functional analysis and operator theory was well understood, it was not until 20 years later that Ruelle employed it with great effect in ergodic theory and dynamical systems, in his application to dynamical zeta functions \cite{ruelleinventiones}. 
This viewpoint was highly influential 
in the work of mathematical physicists
(e.g., the well known study by  Cvitanovi\'c in his monumental online tome \cite{cvitanovic}).

With the advent of modern computers, an important component in  ergodic theory and dynamical systems has been the focus on
\emph{explicit computation} of quantities arising in the context of dynamical systems.
Typically,  these approaches are based on the quantity in question being expressed in terms of an associated Ruelle transfer operator, implicitly assumed to act on some  appropriate space of functions, and then a finite dimensional approximation to the operator is used to reduce this to a finite dimensional matrix problem, to be solved numerically.
In contrast, the  approach we will describe is to  exploit \emph{real-analytic}
properties of the underlying dynamical system by introducing Ruelle transfer
operators with strong spectral properties (in particular nuclearity) which allows us to exploit the earlier circle of ideas initiated by  Grothendieck.

	\subsection{A selection of applications}
By way of an appetiser to this approach, 
we  now list a cross section of actual and potential applications to a number of
different areas in mathematics.  We will elaborate these later, but for the purposes of motivation we list them here.
\begin{enumerate}
\item
In number theory, 
a  recent breakthrough in the understanding of the Markov and Lagrange spectra 
$M, L \subset (0,+\infty)$
 from Diophantine approximation 
has been brought about by the work
of 
Matheus \& Moreira  who were able to estimate the dimension of the difference 
$M\setminus L$ of these spectra, obtaining first a lower bound \cite{matheusmoreira} and then an upper bound \cite{matheusmoreira2}. The methods, in both cases, involved the approximation of the dimension of certain fractal sets which would be amenable to the techniques developed in \cite{jpeffective} for proving rigorous high quality bounds on the dimension
(cf. Example \ref{deleteddigitsexample} (b)). 
\item
In the field of spectral geometry,
there is a strong tradition of computing eigenvalues of the Laplacian, dating back to pioneering work of Hejhal \cite{hej} using classical methods.  However, McMullen's  approximations for 
the lowest eigenvalue for  certain infinite volume hyperbolic manifolds
were based on the dimension 
of the limit set and using this viewpoint could be accurately estimated  \cite{mcmullen3}.  
\item
Within dynamical systems, among the most widely studied numerical quantities are Lyapunov exponents,
 measuring the exponential instability of solutions to various problems. 
There is also interest in estimating Lyapunov exponents
in the theory of random matrix products; 
in particular, in information theory this leads to
the computation of entropy rates for binary symmetric channels, related to examples of hidden Markov chains \cite{holliday}. 
\end{enumerate}

We will develop this last setting  in the following subsection.

\subsection{Illustrative Example: Lyapunov exponents for Bernoulli interval maps}\label{illustrative}
We begin by specifying  a suitable class of hyperbolic transformations.
In particular, we want to assume that the systems we are studying are both real analytic and uniformly hyperbolic,
for example either 
a real analytic expanding map, or
 a real analytic hyperbolic diffeomophism or flow on a locally hyperbolic set.

Before formulating statements in greater generality, let us consider a very specific example of a one dimensional map.
Let 
$X = \mathbb R/\mathbb Z$ be the unit circle and let 
$T: X\to X$ be a piecewise $C^\omega$ Bernoulli  map of the interval which is expanding, i.e., there exists $\lambda > 1$ such that
$|T'(x)| \geq \lambda$ for all $x$.  In this setting it is well known that there is a unique
 $T$-invariant probability measure $\mu$
which is absolutely continuous with respect to Lebesgue measure
(i.e., $\frac{d\mu}{dx} \in L^1(X)$) by the famous Lasota-Yorke theorem (see \cite[\S 5.1]{kh}).

\begin{example}\label{doubling}
	More concretely, suppose $T:\mathbb R/\mathbb Z \to \mathbb R/\mathbb Z$ is defined by
	$T(x) = 2x + \epsilon \sin(2\pi x) \pmod 1$ where $|\epsilon| < \frac{1}{4\pi}$, so  in particular
	$|T'(x)| >  (2 - 2\pi \epsilon) > \frac{3}{2}>1$ for all $x\in \mathbb R/\mathbb Z$.
	\end{example}
	
	  The Lyapunov exponent of the unique absolutely continuous $T$-invariant probability measure $\mu$
	  is defined by
	  \begin{equation}\label{(1.1)}
	  L(\mu) = \int \log |T'(x))| \, d\mu(x)\,,
	  \end{equation}
and by the well known Rohlin identity this also equals the entropy $h(\mu)$.

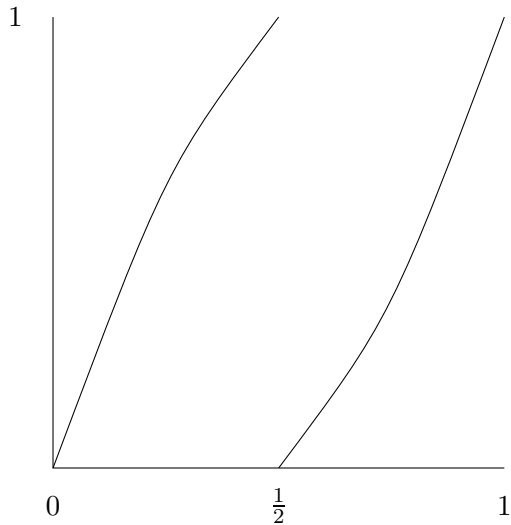
\begin{figure}
\centerline{
\begin{tikzpicture}
\draw (0,6) -- (0,0) -- (6,0);
\draw (0,0) .. controls (1.5,4) .. (3,6);
\draw (3,0) .. controls (4.5,2) .. (6,6);
\node at (0,-0.5) {$0$};
\node at (3,-0.5) {$\frac{1}{2}$};
\node at (6,-0.5) {$1$};
\node at (-0.5,6) {$1$};
\end{tikzpicture}
}
\caption{The graph of an expanding map of $\mathbb R/\mathbb Z$ represented on the unit interval.}
\end{figure}

\medskip
We briefly  summarise the method for estimating $L(\mu)$ in three steps.

\medskip
\noindent
{\bf Step 1 (Complex functions and coefficients)}.
We wish to consider period-$n$ points 
 $T^nx=x$ 
 and then define for $t\in\mathbb{R}$ the  
 coefficients $a_1(t), a_2(t), \ldots$ using the Taylor expansion 
 $$
 \exp \left(
 -\sum_{n=1}^\infty \frac{z^n}{n}\sum_{T^nx=x} 
 \frac{|(T^n)'(x)|^{-t}}{1- ((T^n)'(x))^{-1}} 
 \right) 
 = 1 + \sum_{n=1}^\infty a_n(t) z^n + \cdots .
 $$
 
 \medskip
 \noindent
 {\bf Step 2 (Coefficients and Lyapunov exponents)}.
It can be shown that the Lyapunov exponent $L(\mu)$ given by (\ref{(1.1)}) admits the alternative formulation
\begin{equation*}\label{(1.2)}
 L(\mu)  = 
\frac{ - \sum_{n=1}^\infty a_n'(0)}{\sum_{n=1}^\infty n a_n(0)} ,
\end{equation*}
where both the numerator and denominator are absolutely convergent series.
Truncating these series to give the computable quantity
$$
L_N = \frac{ - \sum_{n=1}^N a_n'(0)}{\sum_{n=1}^N n a_n(0)} \,,
$$
we note that 
$$
L_N \to L\, \hbox{ as } N \to +\infty\,,
$$
and so for a large natural number $N$, the quantity
 $L_N$ is an approximation to $L$.

 \medskip
 \noindent
 {\bf Step 3 (Error bounds)}.
 The quality of the approximation to $L$ given by $L_N$ may be indicated heuristically by comparing how closely $L_N$ and $L_{N-1}$ agree,
though more accurate errors in the approximation can be obtained using:
\begin{enumerate}
\item A certain value $\theta\in (0,1)$, which we refer to as the \emph{contraction ratio}, measuring the extent to which a complex disc $D$ is mapped inside itself by the inverse branches $T_j$ of $T$; 
\item The integrals 
$\beta_k = \frac{1}{r^{2k}} \int_{0}^1
\left|
 \sum_{j} 
 T_j'(x + r e^{2\pi i t})
 (T_j(x + re^{2\pi i t}) )^k
 \right|^2
 dt$, where $1 \leq k \leq L$; and 
 \item The weights $\alpha_k = 
 \sqrt{\sum_{l=k+1}^L \beta_l^2}$ for $N \leq k \leq L$,
\end{enumerate}
 where $N <L$ are suitably chosen.
In particular,  we can then use these values to bound the coefficients $a_n(t)$, with $n > N$,  for which it is impractical to  explicitly compute them with effective error estimates.  This will be explained in greater detail in \S 5.


\bigskip
Now that we have illustrated the general theme using
the specific example of Lyapunov exponents for one dimensional expanding maps,
we can turn to more general dynamical settings, and more general characteristic values.
 In the next section we  
describe the broad context of the results 
and 
consider the pressure function, from which many of the quantities we want to consider can be derived.

\section[Hyperbolic maps]{Hyperbolic maps and the pressure function}
To set the scene, we  begin by introducing two natural classes of discrete dynamical system, then consider the associated pressure function, 
which will be useful in providing the bridge between the dynamics and the various quantities we wish to describe.

 \subsection{Hyperbolic maps}
 Let us begin with the discrete setting.  
 Assume that $T: M \to M$ is either a  smooth 
 expanding map (perhaps on an invariant repeller $Y \subset M$) or Anosov diffeomorphism (see \cite[\S 6.4]{kh}).
 Later we will consider the more general settings of repellers and the natural generalisation to flows.  
However, for the present we will restrict to the discrete cases above and, whenever more convenient, to the case of expanding maps.

\begin{definition}
A partition $X = \cup_{i=1}^k X_i$ 
is called a  {\it Markov partition} for the expanding map $T: X \to X$ if
\begin{enumerate} 
\item
$ \overline{\hbox{\rm int}(X_i)} \cap  \overline{\hbox{\rm int}(X_j)}$ for $i \neq j$;
\item
$X_i = \overline{\hbox{\rm int}(X_i)}$ for $i=1, \cdots, k$; and 
\item each $T(X_i)$ for $i=1, \cdots, k$ is a union of  other elements of the partition.
\end{enumerate} 
\end{definition}


Eventually, we will want to assume  that each of the restrictions 
 $T|X_i$ ($i=1, \cdots, k$) is real analytic, in the sense 
 of having an analytic extension (via charts) to a complex neighbourhood $U_i$.  However, to set up the definitions we only require that it be $C^1$.

In the case of Anosov diffeomorphisms the approach is somewhat similar, except that one uses Markov partitions for invertible maps.  In the case of Anosov flows  one can expect to use Markov Poincar\'e sections to reduce the analysis to the discrete case.

\subsection{Pressure}
The pressure function was introduced into the study of hyperbolic dynamical systems by Ruelle (see e.g.~\cite{ruellebook}).
The importance of pressure stems from the fact that it yields a
unifying concept to describe dynamical and geometric invariants. 
For example, it is well known that various dynamically and geometrically defined fractals (e.g.~limit sets and Julia sets)
have the property that
their Hausdorff dimension is given by solving an associated \emph{pressure equation} (usually known as the Bowen formula).  More generally, there are a host of other dynamical quantities that can be expressed in terms of the pressure function, some of which are listed below in subsection \ref{relating}.


To define pressure $P$, since we are considering hyperbolic maps we have the luxury of expressing this in terms of 
 periodic orbits 
 $T^nx=x$, for $n \geq 1$, and define
$P: C^0\left(\coprod_i X_i\right) \to  \mathbb R $ ($i=1, \cdots, k$)
on the disjoint union of the elements of the Markov Partition.
 
 \begin{definition}
The {\it pressure} of the continuous function $g$ is given by 
 $$
 P(g) := \limsup_{n \to +\infty}
\frac{1}{n}
\log\left( 
 \sum_{T^nx=x} 
 \exp \left( \sum_{j=0}^{n-1} g(T^ix) \right)
 \right)\,,
 $$
 and admits the
 alternative variational definition
$$
P(g) = \sup\left\{ h(\mu) + \int g d\mu \hbox{ : } \mu  \hbox{ is a $T$-invariant 
probability measure}\right\}.
$$
\end{definition}

When $g$ is H\"older continuous, there is a unique probability measure $\mu_g$ realising the above supremum \cite{bowenbook}.

\begin{definition}
The measure $\mu_g$ is called the {\it equilibrium measure} (or {\it Gibbs measure}) for $g$.
\end{definition}

If $g=0$ then $P(0)=h_{\rm top}(f)$ is the topological entropy, and the corresponding equilibrium measure $\mu_g$ is called the measure of maximal entropy (and the Bowen-Margulis measure in the case of Anosov systems).

\begin{example}\label{acim}
	For expanding maps $T:M \to M$, 
	in the special case $g(x)=-\log |\hbox{\rm Jac}(D_xT)|$ then $P(g)=0$, and the corresponding equilibrium measure $\mu_g$ is 
	a $T$-invariant probability measure equivalent to the volume on $M$.
	For Anosov maps $T:M \to M$, 
	in the special case that $g(x)=-\log |\hbox{\rm Jac}(D_xT|E^u)|$, where $E^u$ is the unstable bundle,  then $P(g)=0$  and the corresponding equilibrium measure $\mu_g$ is called the Sinai-Ruelle-Bowen measure
	(or SRB-measure) (see \cite{kh}, \S20.4).
\end{example}

It is this  pressure function that often  helps to relate periodic orbits to the quantities in which we are interested, and which we ultimately want to numerically estimate.
A simple, but important, application  is the following result due to Ruelle
(see \cite[p. 99]{ruellebook}):

\begin{lemma}\label{ruellederivatives}
For any H\"older continuous functions $g_0, g: X \to \mathbb R$ the function 
$$t \mapsto P( g_0+ tg) \in \mathbb R, \hbox{ for $t \in \mathbb R$},$$ is analytic.  Moreover
\begin{enumerate}
\item $\frac{d P(g_0 + tg)}{dt}|_{t=0} = \int g d\mu_{g_0}$, and 
\item 
$\frac{d^2 P(g_0 + tg)}{dt^2} |_{t=0}  = \lim_{n \to +\infty} \frac{1}{n}
\int \left(\sum_{i=0}^{n-1} g(T^ix) \right)^2 d\mu_{g_0}(x)$
provided $\int g d\mu_{g_0}=0$.
\end{enumerate}
\end{lemma}

\smallskip

   The quantity in part 2 of Lemma \ref{ruellederivatives} is often called the {\it variance}.

It is a feature of the method we use that we can obtain fairly explicit  expressions for derivatives of pressure.
In particular, those quantities that can be written in terms of the derivative expressions can therefore, in turn, be written  in terms of periodic points.

\subsection{Relating  pressure to characteristic values}\label{relatingpressure}\label{relating}

We can now consider a number of familiar quantities that we can write in terms of the pressure function.
Below we list a few   simple examples.  Later we shall consider other applications, but for the present these three examples
illustrate well our theme.

\medskip
\noindent
{\bf (I) Lyapunov exponents.}  For an expanding map $T: M \to M$ we can write the 
Lyapunov exponent for the absolutely continuous 
invariant measure $\mu$ by 
$$
L(\mu)= \int \log \|D_xT\| d\mu(x). 
$$   
In the particular case of one dimension this reduces to the situation described in subsection  \ref{illustrative}.

The following follows immediately from part 1 of Lemma \ref{ruellederivatives}.

\begin{lemma}
	If we let $g_0(x)=-\log |\hbox{\rm Jac}(D_xT)|$ and $g = -  \log \|D_xT\|$ then 
	$
	L(\mu) = \frac{d}{dt} e^{P(g_0 + tg)}|_{t=0}.
	$
\end{lemma}

\medskip
\noindent
{\bf (II) Variance.}  For an expanding map $T: M \to M$ we can write the variance  for the absolutely continuous 
$T$-invariant measure $\mu$ 
and a H\"older continuous function $g: X \to \mathbb R$ 
with  $\int g d\mu=0$
defined by
$$
\Sigma(g, \mu) := \lim_{n \to +\infty} \frac{1}{n}
\int \left(\sum_{i=0}^{n-1} g(T^ix) \right)^2 d\mu(x).
$$   
The following follows immediately from part 2 of Lemma \ref{ruellederivatives}.  

\begin{lemma}
We can write 
	$$
\frac{d^2 P(g_0 + tg)}{dt^2}|_{t=0} = \Sigma(g, \mu).
	$$
\end{lemma}

 This plays an important role in the Central Limit Theorem \cite{bowenbook}, i.e.,  for any real numbers $a< b$ we have 
 $$
\lim_{n \to +\infty} \mu\left(\left\{ x \in X \hbox{ : } \frac{1}{\sqrt{n}}\sum_{i=0}^{n-1} g(T^ix) \in [a,b]\right\}\right) = 
\frac{1}{\sqrt{2\pi}} \int_a^b e^{-u^2/2\sigma} du.
 $$

\medskip
\noindent
{\bf (III) Linear response.} 
Let  $T_\lambda: M \to M$ be a smooth family of expanding maps ($-\epsilon < \lambda < \epsilon$)
and let $\mu_{g_\lambda}$ be the associated absolutely continuous measure, arising as the 
Gibbs measure for 
$g_\lambda(x) = 
-\log |\hbox{\rm Jac}(D_xT_\lambda)|
$.
	Then by  part 2 of Lemma \ref{ruellederivatives} we have
	$$\frac{\partial P(g_\lambda + t g)}{\partial t}
	|_{t=0} =   \int g d\mu_{g_\lambda}.$$
Thus differentiating in $\lambda$ formally gives:
$$
\frac{d^2 P(g_\lambda + t g)}{dt\, d\lambda} |_{\lambda=0}=\frac{d}{d\lambda} 
\left(\int g d\mu_{g_\lambda}\right)|_{\lambda=0}.
$$
(A slight subtlety here is that the differentiation of the pressure is easier with respect to the fixed expanding map 
$T_0: M \to M$ and thus it is appropriate to introduce  a family of conjugacies $\mu_\lambda: M \to M$
between $T_0$ and $T_\lambda$ and to consider $\frac{d^2 P(g_0\circ \pi_\lambda + tg\circ \pi_\lambda)}{dt d\lambda} |_{\lambda=0}$).



The above list does not exhaust the possible quantities that can be derived from the pressure, but gives a selection 
we hope illustrates our general approach.
 
\bigskip
In the next section, we  will introduce a standard tool, the  \emph{transfer operator}, which allows us to analyse the pressure, and thus its many derivative properties, using basic ideas from linear operator theory.  
We will also describe the connection with a family of complex functions called determinants.

\section[Transfer operators]{Transfer operators and determinants}

A central object in thermodynamic formalism is the \emph{transfer
operator}, from which important dynamical and geometric invariants
such as entropy, Lyapunov exponents, invariant measures, and Hausdorff
dimension can be obtained.


Let us now  restrict (for the present) to the case of expanding maps $T: X \to X$.
The analyticity of the pressure, as well as other properties  including the 
proof of Lemma \ref{ruellederivatives}, depend on the use of {\it transfer operators}.  
Eventually, we will want to consider operators acting on spaces of analytic functions, but for the purposes of defining them it suffices for the present to consider the Banach space of $C^1$ function $C^1(X, \mathbb C)$ with the norm 
$\|f\| = \|f\|_\infty + \|Df\|_\infty$. 
The operators are  then defined as follows:

\begin{definition}  
If $T:X \to X$ is a $C^1$ expanding map,
and $g: X \to \mathbb R$ is $C^1$,
 then we define the {\it transfer operator}
 $\mathcal L_{g}: C^1(X, \mathbb C) \to C^1(X, \mathbb C)$ by 
$$
\mathcal L_g w(x) = \sum_{Ty = x} e^{g(y)} w(y) \,,
$$
the summation being over the inverse images $y$ of the point $x \in X$.
\end{definition}

 This operator preserves various function spaces and exhibits
certain positivity properties which ensure that a Perron-Frobenius
type theorem holds: 
when acting on $C^1(X, \mathbb C)$,
$\mathcal L_g$ has a leading eigenvalue which is simple, positive, and isolated.
Moreover, the connection with the pressure comes from the following basic result due to Ruelle \cite{ruellebook}
(see also  \cite{bowenbook}).
\begin{lemma}
	The spectral radius of $\mathcal L_g$ is 
$e^{P(g)}$.	
	In particular,  
$e^{P(g)}$ is a maximal isolated eigenvalue for 
$\mathcal L_g$.
\end{lemma}

In particular, the differentiability 
(indeed analyticity) of the pressure follows by standard perturbation theory and the expression for the derivatives in the lemma follow by explicit manipulations.

\begin{example}
In the special case that $g(x)=-\log |\hbox{\rm Jac}(D_xT)|$ then $P(g)=0$  and 
$\mathcal L_{g}$ is known as the {\sl Ruelle-Perron-Frobenius operator}.
The eigenmeasure $m_g = \mathcal L_g^*m_g$ is normalised Lebesgue measure,
and the equilibrium measure $\mu_g = h_g m_g$ is the unique
$T$-invariant measure absolutely continuous
with respect to Lebesgue measure, where  $h_g = \mathcal L_g h_g $ is the maximal eigenfunction \cite{bowenbook}.
\end{example}

Thus far we have been following a very traditional approach.  However, now we introduce an extra  ingredient.

\subsection{Determinants and their coefficients}
Let $T: X \to X$ be a $C^1$ expanding map.
 For any continuous function $G: X \to \mathbb R$ 
and each period-$n$ point $T^nx=x$, $n \geq 1$, we can associate the weight 
$$ G^n(x):= \sum_{i=0}^{n-1} G(T^ix) \in \mathbb R.$$
Later we will want to assume that $T$ and $G$ are real analytic, but for the purposes of introducing the determinant we need 
only assume these weaker hyotheses.
It is convenient to package up the information from individual periodic points into a single (generating) complex  function.

\begin{definition}
Given a continuous  function $G: X \to \mathbb R$ we can formally define  a function of the single complex variable
by:
$$
D(z)  =  D_{G,T}(z) = \exp\left(- \sum_{n=1}^\infty 
\frac{z^n}{n} \sum_{T^nx=x} 
\frac{\exp\left({\sum_{i=0}^{n-1}G(T^ix)}\right)}{
	\det(I- [D(T^n)(x)]^{-1})}
\right), \quad z\in \mathbb C\,.
$$
\end{definition}
\bigskip

The radius of convergence of the infinite series $D(z)$ is related to the pressure.  More precisely, 
we can see that this converges to an analytic function provided the series converges, i.e., 
$|z| e^{P(G)} < 1
$ 
where $$e^{P(G)} = \lim_{n\to +\infty} \left|\sum_{T
	^nx=x}  \frac{\exp\left(\sum_{i=0}^{n-1}G(f^ix)\right)}{\det(I- [D(T^n)(x)]^{-1})}\right|^{1/n}
\left(=  \lim_{n\to +\infty}
\left|\sum_{T
	^nx=x}  \exp\left(\sum_{i=0}^{n-1}G(f^ix)\right)\right|^{1/n}
\right).
$$
In particular, writing $D(z)$ as a power series  
$$
D(z) = 1 + \sum_{n=1}^\infty a_n z^n,\eqno(3.1)
$$
with coefficients $a_n = a_n(T,G)$  depending on $T$ and $G$,  we see it  has  radius of convergence
at least $e^{-P(G)}$.

\begin{example}[Expanding maps of the interval]
In the particular case of an expanding map $T: X \to X$ of the interval $X$, 
given a continuous  function $G: X \to \mathbb R$,  the function $D(z)$ takes the simpler form:
$$
D(z)  =   \exp\left(- \sum_{n=1}^\infty 
\frac{z^n}{n} \sum_{T^nx=x} 
\frac{\exp\left({\sum_{i=0}^{n-1}G(T^ix)}\right)}{
 1- 1/(T^n)'(x)}
\right), \quad z\in \mathbb C\,.
$$
\end{example}

This naturally leads to asking about the meromorphic extension of $D(z)$.
To proceed further, we need to assume more regularity 
on the function $G$.  
This  brings us to the following important result of Ruelle \cite{ruelleinventiones}.

\begin{lemma}\label{ruelle-domain}
	If $T$ is real analytic then
	\begin{enumerate}
		\item
 $D_{G,T}(z)$ is analytic in all of  $\mathbb C$.
	\item
	The value $z= e^{-P(G)}$ is a simple zero for $D_{G,T}(z)$ in this extension.
	\end{enumerate}
\end{lemma}

In particular, we see from 	part 1 of Lemma \ref{ruelle-domain}  that we can improve the result on the radius of convergence of the power series to     $\lim_{n\to +\infty} |a_n|^{1/n} = 0$, i.e., for any $0 < \theta < 1$ 
there exists $C>0$ such that $|a_n| \leq C \theta^n$.
In fact, the original proof of 
Lemma \ref{ruelle-domain} due to Ruelle, and inspired by work of Grothendieck, 
works in this way by giving estimates on the coefficients $a_n$.  We will later describe  quite precise  bounds on the coefficients $a_n$ which establishes part 1. We will return to this point in the next subsection.



 \begin{rem}
 If we assume that $T$ and $G$ are $C^\infty$ then we would still have that $\lim_{n\to +\infty} |a_n|^{1/n} = 0$.  
 However, as we shall see,  in the analytic case we have more effective estimates on $|a_n|$.
 \end{rem}
 
 \subsection{Pressure and the characteristic quantities}
 In order to relate $D(z)$ back to the pressure, and thus the various dynamical quantities, we need to make different choices for $G$.
 More precisely, we  can consider  the special case of the function 
 $G = g_0 + t g$ where $g_0, g: X \to \mathbb C$ and $t\in \mathbb R$.   This leads to the following particular  case of the previous definition.

 \begin{definition}\label{defn}
 	We  formally define the {\it determinant} for $g_0$, $g_1$ to be the bi-complex function
 $$
 d_{g_0,g}(z, t):= 
 \exp \left(
 -\sum_{n=1}^\infty \frac{z^n}{n} \sum_{T^nx=x} 
\frac{ \exp \left( \sum_{j=0}^{n-1} (g_0 + tg)(T^ix) \right)
}{
1 - 1/(T^n)'(x)
}
 \right),
 $$
 where $z,t \in \mathbb C$.
 \end{definition}

The relationship between the determinant and the pressure 
in Lemma \ref{ruelle-domain}
now implies the following.

\begin{cor}\label{cor}
	Assume that   $g_0$, $g_1$ are $C^\omega$.
	\begin{enumerate}
		\item
		The function 
 $d_{g_0,g}(z, t)$ has a bi-analytic  extension to  all of  $\mathbb C^2$.
	\item
	The value $z= e^{-P(g_0 + tg)}$ occurs as a simple zero for $z \mapsto d_{g_0,g}(z, t)$.
	\end{enumerate}
\end{cor}

We can make the further  simplifying assumption that $P(g_0)=0$, where we can replace $g_0$ by $g_0 - P(g_0)$ if necessary.  In particular, the first zero for $z \mapsto d_{g_0,g}(z, 0)$ (where $t=0$) occurs at $z=e^{-P(g_0)}$.  
We can use part 2 of Corollary \ref{cor}, and the implicit function theorem, to write the derivative of the pressure as
$$
\frac{d P(g_0 + tg)}{dt}|_{t=0} = 
\frac{\partial d_{g_0,g}(1, t)}{\partial t}|_{t=0}/
\frac{\partial d_{g_0,g}(z, 0)}{\partial z}|_{z=1}
$$
in terms of partial derivatives of the determinant.
 Furthermore, in light of part 1 of Corollary \ref{cor}, for each $t \in \mathbb R$  we can formally expand 
 $$
  d_{g_0,g}(z, t) = 1 + \sum_{n=1}^\infty a_n(t) z^n.
 $$
 It is these values which we need to relate to the quantities in which we are interested, and which we ultimately want to numerically estimate.
 
 This  gives us the simplest definition of the 
 numbers $(a_n(t))_{n=1}^\infty$, 
 although we can explicitly expand these in terms of periodic points too.  This reveals the following simple, but crucial, fact.

\begin{lemma}
For each $n \in \mathbb N$, we can express the value $a_n(t)$ in terms of the periodic points of period at most $n$.
\end{lemma}
 
 In particular, this ensures that the coefficients 
$a_n(t)$ are relatively easy to estimate.  
 Moreover, we can rapidly approximate 
 $
 d_{g_0,g}(z, t) = 1 + \sum_{n=1}^\infty a_n(t) z^n
 $
 by the truncated series
 $$
 1 + \sum_{n=1}^N a_n(t) z^n
 $$
 for $N$ moderately large.
  The rapidity of the approximation is justified by the following result.

\begin{cor}\label{super}
 The coefficients $a_n=a_n(t)$ tend to zero at a super-exponential rate.  \end{cor}
 
When $X$ is one dimensional then there exists $\theta \in(0, 1)$ such that 
$|a_n| = O(\theta^{n^2})$ as $n\to\infty$.  
We will give very explicit estimates for the implied constants in the $O(\cdot)$ term.

\begin{rem}
When  $X$ is $d$-dimensional (with $d \geq 2$) then there exists $0 < \theta < 1$ such that 
$|a_n| = O(\theta^{n^{1+1/d}})$ as $n\to\infty$.  
\end{rem}

\begin{rem}
More generally, we can assume we have a family of real analytic 
functions $g_1, \cdots, g_m: X \to \mathbb C$ and 
 $$
 d(z, t):= 
 \exp \left(
 -\sum_{n=1}^\infty z^n \sum_{T^nx=x} 
 \exp \left( \sum_{j=0}^{n-1} (g_0 + t_1g_1 + \cdots +  t_mg_m)(T^jx) \right)
 \right)
 $$
 for $z \in \mathbb C$ and $t_1,\ldots,t_m \in \mathbb R$.
\end{rem}  

To summarize, we now have a  
method of approaching the pressure function which might be considered to have a simple  analogy  to that  of estimating the largest eigenvalue of a matrix  by computing the characteristic polynomial,  which  is a complex function whose zeros give the eigenvalues.  This simple viewpoint  relating  the determinant and transfer operators ultimately leads to a surprisingly efficient method of computing pressure.
 Furthermore, we can accurately estimate the error terms (see \S 5).

Having related the various  quantities which we want to estimate to the zero(s) of the determinant $d(z,t)$ (often via the pressure function, cf. subsection 2.3) we need to  answer three key questions:
 {\it  How can we use this formulation to get numerical estimates?
Why does this approach lead to superior approximation estimates?
How can we estimate the quantities with validated rigour?}

We can now turn to the practical problem of writing explicit expressions for the approximations to quantities we described in subsection  \ref{relating}, in terms of the first $N$ coefficients in the expansion of the determinant.

  \subsection{Dynamical quantities and coefficients}\label{quantities}
We have established  (in subsection 2.3) that many  of the quantities that we want to estimate can be expressed in terms of pressure, and its derivatives, which in turn can be written in terms of the determinant and its derivatives.  
Since the determinant has a power series expansion it is a  straightforward, but useful,  exercise  to write these expressions explicitly in terms of the derivatives of the coefficients.    More precisely,  let us write
$$
A = \sum_{n=1}^\infty n a_n(0),\quad
B = \sum_{n=1}^\infty n(n-1) a_n(0),\quad
C = \sum_{n=1}^\infty a_n'(0),\quad
D = \sum_{n=1}^\infty n a_n'(0),\quad
E = \sum_{n=1}^\infty  a_n^{\prime \prime}(0)
$$
and the associated finite sums
$$
A_N = \sum_{n=1}^N n a_n(0),\ \
B_N = \sum_{n=1}^N n(n-1) a_n(0),\ \
C_N = \sum_{n=1}^N a_n'(0),\ \
D_N = \sum_{n=1}^N n a_n'(0),\ \
E_N = \sum_{n=1}^N  a_n^{\prime \prime}(0)
$$
for $N \geq 1$.  

A recurrent theme in our discussions   is that 
we want to express the dynamical quantities in terms of 
$A$, $B$, $C$, $D$, $E$, etc.,  and then approximate these expressions by using instead the more computationally tractable quantities 
$A_N$, $B_N$, $C_N$, $D_N$, $E_N$.

\begin{rem}
In practical applications, even for quite simple examples,  
we might currently only expect to compute these values up to $N=25$, say, in a reasonable time frame.
However, with further technological advances, one might expect that this value can be improved.
\end{rem}





To illustrate this principle,  we can  now reformulate the three key quantities described in subsection \ref{relatingpressure}, and their approximations, in terms of these series and summations.
We list these below.

\medskip
\noindent
{\bf (I) Lyapunov exponents.}
We can write the Lyapunov exponent for $\mu$ as
$$L(\mu)  = -\frac{\sum_{n=1}^\infty a_n'(0)}{\sum_{n=1}^\infty n a_n(0)}
= -\frac{C}{A}\,,
$$
and in particular the Lyapunov exponent can be approximated by the computable quantities
$$
-\frac{C_N}{A_N}, \quad N \geq 1\,.
$$

\medskip
\noindent
{\bf (II) Variance.}
The variance is given by
$$
\Sigma^2 = \left(\frac{C}{A}\right)^2 
+ \frac{1}{A}
\left(B \left(\frac{C}{A}\right)^2 - 2 D B \left(\frac{C}{A}\right) + E\right) \,,
$$
and in particular we can approximate the variance by
$$
 \left(\frac{C_N}{A_N}\right)^2 
+ \frac{1}{A_N}
\left(B_N \left(\frac{C_N}{A_N}\right)^2 - 2 D_N B_N \left(\frac{C_N}{A_N}\right) + E_N\right), 
\quad N \geq 1\,.
$$


\medskip
\noindent
{\bf (III) Linear response.}
We can write
$$
\int g d\mu_T  = -
\frac{C}{A} \,,
$$
and in particular we can approximate the integral by 
$$
-\frac{C_N}{A_N}, \quad N \geq 1.
$$
Finally, by replacing $T$ by $T_\lambda$, differentiating both sides in $\lambda$ and we get an expression for the linear response in terms of the (derivatives of the) coefficients $a_n$.



 
\section{Rates of mixing and dimension}\label{mixing}
In this section we want to 
make a slight detour to 
introduce another two important quantities which,
although not quite fitting into the same framework  described in the previous section,
 can also be approximated using the determinant. 


First we consider the rate(s) of mixing, which can be studied via the zeros of the determinant.


\subsection{Rates of mixing}\label{rates}
Let  $X$ be $d$-dimensional and let $T: X \to X$ be  a $C^\omega$ conformal  expanding map.  More precisely, we can write the derivative
$DT(x) = w(x) \Theta(x)$ where $w: X \to \mathbb R$ and $\Theta: X \to SO(d)$.   In the particular case that $d=1$ then the one dimensional map 
$T$  is automatically conformal.

Let $g_0: X \to \mathbb R$ be  real analytic and 
let $\mu = \mu_{g_0}$ be the equilibrium-Gibbs measure 
associated to $g_0$.

\begin{example}As we observed in Example \ref{acim}, when $g_0(x) = -\log |\hbox{\rm Jac}(T)(x)|$ the associated measure $\mu_{g_0}$ is the unique absolutely continuous invariant probability measure. 
\end{example}

The following is an important object in ergodic theory.

\begin{definition}
 Given a real analytic function
  $g$, we can consider the {\it correlation function} defined by
  $$c(n)= 
  \int g \circ T^n g d\mu - \left(\int g d\mu\right)^2, \quad n \geq 1.
  $$ 
  \end{definition}
  
  Since $T:(X, \mu) \to (X, \mu)$ is mixing we know that $c(n) \to 0$,  as $n \to +\infty$.
    The (exponential) rate of mixing is given by the smallest  value $0 < \lambda_1 < 1$
such that $c(n)= O(\lambda_1^n)$ for all such $g$ and $n \geq 1$.
%
Since $\lambda_1$  corresponds to the modulus of the second eigenvalue  of the transfer operator, 
the connection to the determinant comes through the following  simple lemma.

\begin{lemma}\label{rho}
The rate of mixing $0 < \lambda_1 < 1$ is the reciprocal of the modulus $\rho>1$ of the second smallest (in modulus) zero of 
$d_{g_0,g}(z,0)$, i.e.~$\lambda_1 = 1/\rho$.
\end{lemma}

In particular, the value $\rho$ in Lemma \ref{rho} is a zero of 
the (real valued) series
$$
d_{g_0,g}(z,0)
= 1 + \sum_{n=1}^\infty z^n a_n(0)
$$
and so in order to get rigorous bounds on $\rho$ we can use the intermediate value theorem.
More precisely, given $\epsilon_1, \epsilon_2>0$ and
 $N \in \mathbb N$
we  look  for  bounds $\alpha_N < \rho < \beta_N$ by choosing 
$\alpha_N, \beta_N$ such that 
$$
1 + \sum_{n=1}^N \alpha_N^n a_n(0) \leq  -\epsilon_1
\hbox{ and }
1 + \sum_{n=1}^N \beta_N^n a_n(0) \geq  \epsilon_2
$$
and
$$
\epsilon_1 > \left|\sum_{n=N+1}^\infty \alpha_N^n a_n(0)\right|
\hbox{ and }
\epsilon_2 > \left|\sum_{n=N+1}^\infty \beta_N^n a_n(0)\right|.
$$

\noindent 
Thus finding good bounds $\alpha_N < \rho < \beta_N$ reduces to:
\begin{enumerate}
\item
getting 
good estimates on the coefficients $a_i(0)$ 
($i=1, \cdots, N$); and 
\item 
finding  effective bounds on the tails
$\sum_{n=N+1}^\infty a_n(0)$.
\end{enumerate}
The first is a basic problem in practical  computing.  The second is a more challenging mathematical problem.

We can illustrate this with a simple example.

\begin{example}[Lanford map, cf.~\cite{jpv,lanford}]\label{lanford}
  Let $T: [0,1] \to [0,1]$ be defined by 
$$
T(x) = 2x + \frac{1}{2} x(1-x) \pmod 1\,.
$$
Let $\mu$ denote the unique absolutely continuous $T$-invariant probability measure.
We can now evaluate dynamical quantities such as the rate of mixing
by approximating the infinite series for the determinant by finite truncations.
For example, using periodic points up to period 
$N= 16$  we might estimate the first ten zeros (i.e.~the ten zeros of smallest modulus)
of  the determinant  by
locating the first ten zeros of the degree-16 polynomial truncation
(see Table 1).

\begin{table}[h!]
  \begin{center}
    \begin{tabular}{| c | l |} 
    \hline
    $n$ & zero $z$\\ 
    \hline
    $1$ & $1.0000000000000000000000033711203720152 $\\
$2$& $1.72986531066431681927894069519181629 $\\
$3$ & $2.6922698183465737455729975178528581 $\\
$4$&$4.1132466756759777783645672672979956 $\\
$5$& $6.2454853205721291176033177124804291 $\\
$6$ & 9.4538916717397326544473431332123370 \\
$7$ & $ 14.282734336458524434000313080802510 $ \\
$8$  & $21.549229994532327179757991408669084$ \\
$9$   & $32.516266102701803490675636630907193$ \\
$10$ & $47.82910484702218758446773289813753 $\\
\hline
 \end{tabular}
  \end{center}
  \caption{Estimates on the  first $10$ zeros of the determinant}
\end{table}

More accurate estimates on the  second eigenvalue $\lambda_1$ of the transfer operator,
i.e.~the exponential rate of mixing, 
are obtained by approximating the reciprocal of the second zero of the determinant, given in Table 2,
using degree-$N$ truncations with $N$ ranging up to $N=25$.

\begin{table}[h!]
    \begin{tabular}{| c | l |} 
    \hline
    $N$ & second eigenvalue $\lambda_1$ estimate\\ 
    \hline
12&
0.5780796887515271422742765368788953299348846128812023109203951947004787498004165\\
13&
0.5780796885356306834127405345836355663641109763750019611087170244976104563627485\\
14&
0.5780796885371288506764371131157188309769151998850254045247866596035386808066373\\
15&
0.5780796885371219470570630291328371225537224787114789966418506438634692131905786\\
16&
0.5780796885371219681960432055118626393344991913606205477442507113706445878179303\\
17&
0.5780796885371219681530107872274433995896003010891980049121721575572541941602200\\
18&
0.5780796885371219681530690475964044549434630264578745046610737538545621059499499\\
19&
0.5780796885371219681530689951228630661230046404665596121577924787108887084624467\\
20&
0.5780796885371219681530689951543111607290086789469044407358342311002102577614591\\
20&
0.5780796885371219681530689951543111607290086789469044407358342311002102577614591\\
21&
0.5780796885371219681530689951542986173994750355886730998629996708085266331162364\\
22&
0.5780796885371219681530689951542986207295902353572376170142239221705813115693187\\
23&
0.5780796885371219681530689951542986207290016813780055122529428636713090405312973\\
24&
0.5780796885371219681530689951542986207290017506309910801987018396260494147990458\\
25&
0.5780796885371219681530689951542986207290017506255654011865278736341388970140860\\
\hline
\end{tabular}
\caption{Estimates on the second eigenvalue $\lambda_1$ coming from the reciprocals of zeros for $1 + \sum_{n=1}^N z^n a_n(0)$, for $12\leq N \leq 25$.}
\end{table}

\end{example}

We will address the rigorous bounds on the error later.

\begin{rem}
The rate of mixing has a more subtle generalisation of the following form.
The speed of mixing is controlled by  a sequence of complex numbers $\lambda_j \to 0$ 
ordered to be decreasing in modulus, and polynomials $P_j$.   More precisely, 
for any $\delta > 0$ there exist $M$ and $c_j = c_j(g)$ ($j=1, \cdots, M$) such that 
$$
c(n) = \sum_{j=1}^M c_j P_j(n)\lambda_j^n 
$$
where $P(n)=1$  in the case of a single zero of the determinant.   
  Since we can rewrite
  $$
  c(n) = \int g \left( \mathcal L_{g_0}^n  g\right)  d\mu - \left( \int g  d\mu\right)^2
  $$
  the values  $\lambda_j$ are actually the 
   {\it other} non-zero eigenvalues
  of the transfer operator.   Furthermore,
  $$
P_j(n) =   
\begin{cases}
 1 & \hbox{ if  } \lambda_j \hbox{ has multiplicity }m  = 1\\
 \hbox{ a polynomial of degree at most $n$} & \hbox{ if  } \lambda_j \hbox{ has multiplicity $>1$}.\\
\end{cases}
$$
Furthermore, one can identify the eigenvalues in terms   of the {\it other} zeros $\rho_j = 1/\lambda_j$ for 
 $z \mapsto d_{f,z}(z,1)$.  
 
 \noindent
 All this said, in the case of the Lanford map the first 10 zeros of the determinant are simple.
 
 \end{rem}

\subsection{Dimension of repellers}\label{repellerssubsection}
Let  $T: X \to X$ be  a $C^\omega$ conformal  expanding map.

\begin{definition}
A closed invariant set  $Y \subset X$ is called a repeller if there is an open set $U$ satisfying $Y \subset U \subset X$ such that 
$Y = \cap_{n=1}^\infty T^{-n}U$.  
\end{definition}

We refer the reader to the book of Falconer for the definition and basic properties of Hausdorff dimension \cite{falconer}.  
In the present context we have a convenient dynamical formulation.
Let $g_0(x) = -\log |\hbox{Jac}(T)(x)|$ and then consider the restriction $T: Y \to Y$.
The following standard result relates the dimension $ \hbox{\rm dim}_H(X)$ to the pressure function for the transformation 
$T: Y \to Y$ and function  $sg_0: Y \to \mathbb R$, for $s > 0$ (see \cite{ruellerepellers}).

\begin{lemma}[Bowen-Ruelle]\label{br}
We can characterise the 
Hausdorff dimension of the limit set by 
$s = \hbox{\rm dim}_H(Y)$ such that
 $P(sg_0) = 0$.
\end{lemma}

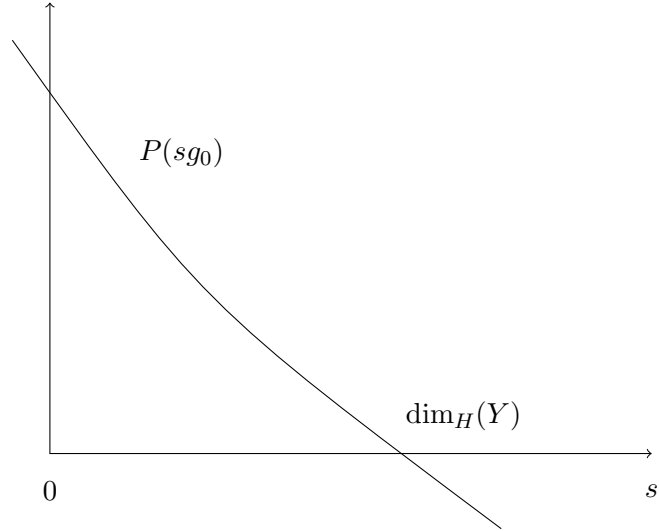
\begin{figure}
\centerline{
\begin{tikzpicture}
\draw[<->] (0,6) -- (0,0) -- (8,0);
\draw (-0.5,5.5) .. controls (2,2) .. (6,-1);
\node at (0,-0.5) {$0$};
\node at (5.5,0.5) {$\dim_H(Y)$};
\node at (1.75,4.0) {$P(sg_0)$};
\node at (8,-0.5) {$s$};
\end{tikzpicture}
}
\caption{Graph of a pressure curve $s\mapsto P(sg_0)$}
\end{figure}

The following examples fit into the setting of Lemma \ref{br}.

\begin{example}
If $T:J\to J$
is a hyperbolic rational map acting on its
Julia set $J$, and $f=-s\log|T'|$, 
then the Hausdorff
dimension $\dim_H(J)$ is given 
by the unique zero $s$ of $P(-s\log|T'|)$ by Lemma \ref{br}.
\end{example}

It is a result of Ruelle \cite{ruellerepellers}
that the  dimension of a hyperbolic  Julia set $J$ varies analytically with $f$. In the interests of clarity, 
for the discussion below we shall restrict attention to the 
specific example of the quadratic map $f_c(z) = z^2+c$ where $c$ is in the main cardioid of the Mandelbrot set.

\begin{thm}[Ruelle]\label{ruelleasymp}
Let $c$ be in the main cardioid of the Mandelbrot set. Then the mapping $c \mapsto \dim_H(J_c)$ is real analytic. 
Moreover,  if  $|c|$ is sufficiently small then there is an asymptotic expansion 
 $$
 \dim_H(J_c) = 1 + \frac{|c|^2}{2\log 2} + O(|c|^3)\,.
 $$

\end{thm}

We briefly sketch the  argument for the  analyticity below.  
Although not central to our survey, for completeness we also include a brief account of the asymptotic expansion in Appendix B.

To make use of Lemma 4.7 we first need that there exists  a holomorphic function $\Phi_c$ which conjugates
 $f_c$ to $f_0$, i.e. $f_0 \circ \Phi_c = \Phi_c \circ f_c$ (see  \cite{Zin00}) and if $z$ has $|z| > 1$ then $\Psi_0(z) = z$.
Moreover, this  function is holomorphic in $c$ for $c$ in the main cardioid 
$\mathcal{C} = \{w (1-w) \hbox{ : } |w| < \frac{1}{2} \}$
of the Mandelbrot set, as is  $\Psi_c = \Phi_c^{-1}$.
Moreover, for $c$ inside this main cardioid:

\begin{enumerate}
\item For any $z$ with $|z|>1$, the map $c \mapsto \Psi_c(z)$ is holomorphic in $c$.
\item For $z$ in the unit circle $J_0$, the map $z \mapsto \log 2|\Psi_c(z)|$ is H\"older continuous.
\item For any $c \in \mathcal{C}$, the map $z \mapsto \Psi_c(z)$ is injective on $\{z \hbox{ : } |z| > 1\}$.
\end{enumerate}
The function $\Psi_c$ extends to a holomorphic map on $\{z \hbox{ : } |z| \geq 1\}$, with $\Psi_c|_{J_0}: J_0 \to J_c$
and  for $z$ on the unit circle, $z \mapsto \log2|\Psi_c(z)|$ is H\"older continuous.
By Lemma 4.7, 
\begin{equation*}
P(-t\log(2|\Psi_c(z)|)|_{J_0}) = 0.
\end{equation*}
 So $(c,t) \mapsto P(-t\log2|\Psi_c|)$ is  analytic in a neighbourhood of $(0,1)$.
We then have:
\begin{equation*}
\frac{\partial}{\partial t} (P(-t\log 2|\Psi_{c=0}(z)|)) |_{t=0}= \int_{J_0} -\log 2|\Psi_{c=0}(z)| \, dz \neq 0.
\end{equation*}

\begin{rem}
Theorem \ref{ruelleasymp} can be generalised to the context of hyperbolic Julia sets of rational maps (see \cite{ruellerepellers}).
\end{rem}

\begin{rem}
The real analyticity of Theorem \ref{ruelleasymp} may fail outside of the main cardioid of the Mandelbrot set.
 For example, when  $c = \frac{1}{4}$ then $ f_c(\frac{1}{2})=  \frac{1}{2}$ and  $|f'(1/2)| = 1$,  so $J_{c}$ is no longer hyperbolic, and in fact $c \mapsto \dim_H(J_c)$ is not even continuous  at $\frac{1}{4}$
 (see \cite{dsz}).
\end{rem}

The same ideas apply to the Hausdorff dimension of limit sets of
certain Fuchsian  groups (see, for example,  \cite{bowenpublihes}).

 Lemma \ref{br} allows us to define the dimension $\hbox{\rm dim}_H(Y)$ implicitly in terms of the  determinant
 (defined in terms of periodic points for $T: Y \to Y$, i.e., those contained in $Y$).  
In particular,  setting $g_0=0$, $g = -\log |T'|$  and $z=1$ in 
Definition \ref{defn}
we have  $s = \hbox{\rm dim}_H(X)$  satisfies 
 $d_{0, g}(1, s) = 0$.  
Equivalently, the 
Hausdorff dimension of the limit set is a solution 
$s = \hbox{\rm dim}_H(X)$ to  the absolutely convergent series
 $$1 + \sum_{n=0}^\infty a_n(s) = 0.$$
As before, when studying the rate of mixing  in the previous subsection,  in practice we can 
use the intermediate value theorem to 
get effective bounds $\alpha_N < \rho < \beta_N$ by choosing 
$\alpha_N, \beta_N$ such that 
$$
1 + \sum_{n=1}^N a_n(\alpha_N) \geq  \epsilon_1
\hbox{ and }
1 + \sum_{n=1}^N a_n(\beta_N) \leq - \epsilon_2
$$
where 
$$
\epsilon_1 > \sum_{n=N+1}^\infty a_n(\alpha_1)
\hbox{ and }
\epsilon_2 > \sum_{n=N+1}^\infty a_n(\alpha_2)
$$
Thus, as before, finding good bounds comes down to:
\begin{enumerate}
\item
getting 
good estimates on the coefficients $a_i(s)$ 
($i=1, \cdots, N$); and 
\item 
finding  effective bounds on the tail
$\sum_{n=N+1}^\infty a_n(s)$.
\end{enumerate}

As we mentioned before, the  first is a basic problem in computing and the second is a more challenging mathematical problem
which we will address in \S 5.

\smallskip
We now turn to a class of deceptively simple examples.

\begin{example}[Continued fractions and deleted digits] \label{deleteddigitsexample}
Consider a finite set $F \subset \mathbb N$ and the set $E_F \subset [0,1]$ given by
$$
E_F = \{x = [x_1, x_2, \cdots ] \hbox{ : }  x_1, x_2, \cdots \in F\},
$$
i.e., the Cantor set of points whose continued fraction expansion has all coefficients lying in $F$.
This can be viewed as a repeller for the map $T: E_F \to E_F$ defined by $T(x) = \frac{1}{x} \pmod 1$.

\medskip
\noindent 
(a)
In the case $F = \{1,2\}$ the set $E_F$ is usually denoted  $E_2$ and is of historical interest,
with its Hausdorff dimension  studied by Good
\cite{good} as far back as 1941, after even earlier work of Jarnik \cite{jarnik}.  Using our algorithm we were able to compute
 the dimension accurately, and rigorously,
 to over 100 decimal places (see \cite{jpeffective}),
improving on the  previous best rigorous estimate due to Falk \& Nussbaum 
\cite{falknussbaum}, using a subtler variant of Ulam's method.

\medskip
\noindent 
(b)
 In the case $F = \{1,2,3,4,5\}$, the  dimension $\dim_H (E_{ \{1,2,3,4,5\}})$ 
 appears as a crucial ingredient in the work of Huang \cite{huang}, 
 refining the work of Bourgain \& Kontorovich \cite{bourgainkontorovich} on a density one solution to the Zaremba Conjecture.
  Here we were able to use the algorithm to compute the dimension accurately, and rigorously, to 8 decimal places
  \cite{jpzaremba}.

\end{example}

\begin{rem}
The Cantor sets above are also special cases of  limit sets of one dimensional  iterated function schemes.
More precisely,  we can consider a (finite) set of $C^2$ contractions $\psi_i:I \to I$ ($i=1, \cdots, N$) on the unit interval $I$ (with $\max_{1\leq i \leq N}\sup_{x\in I} |\psi_i'(z)| < 1$).   The associated limit set $X \subset I$
is the smallest non-empty closed set satisfying $X = \cup_{i=1}^N \psi_i(X)$.  
The Hausdorff dimension $\dim_H(X)$ can then be expressed in terms of the pressure function and the corresponding version of  Lemma \ref{br}.
 In the particular case that 
$\psi_i(x)= \frac{1}{x+i} $, for $i \in F$,  the limit set  recovers the Cantor sets  $E_F$.   
\end{rem}
\begin{rem}
There are natural analogues of continued fractions for which the contractions are $\psi_b(z) = \frac{1}{z+b}$, where 
$b \in \mathcal B \subset \{m + i n \hbox{ : } m \in \mathbb N, n \in \mathbb Z\}$ \cite{mu}.   These map the domain 
$D = \{z\in \mathbb C \hbox{ : } |z-\frac{1}{2}| \leq \frac{1}{2}\}$ inside itself.  We denote by $X_{\mathcal B}$ the associated limit set, i.e., the smallest non-empty closed subset of $D$ such that 
$\cup_{b \in \mathcal B} \psi_b(X_{\mathcal B}) = X_{\mathcal B}$.  The dimension of this set was  studied in \cite{mu} and in \cite{falknussbaum1} the bounds $ 1.85574 \leq \dim_H(X_{\mathcal B}) \leq 1.85590$ were established.
\end{rem}

\begin{rem}
These examples also help to highlight some of the limitations of the periodic point method for estimating values.
Most of the examples we have studied have been one dimensional and involved uniformly expanding maps (or uniformly contracting iterated function schemes) with only finitely many branches.   If we consider analytic maps in higher dimensions then the method still applies, although it is less efficient as the dimension grows.
However, if we consider the case of infinitely many branches then we have the  complication 
that there can be infinitely many periodic points of any given period and this genuinely makes the method less applicable.
\end{rem}

\section{Rigorous error bounds}
We now come to one of the most interesting and challenging aspects of the estimation problem:  {\it finding rigorous upper bounds for  the errors in the approximations.}

In the interests of clarity, and notational simplicity,  we shall 
explain the ideas in the particular case of estimating the Hausdorff dimension of limit sets (corresponding to conformal iterated function schemes).  The more general settings require variants of this basic approach.

In particular, we  need to provide an estimate on the error when we truncate the series
which comes from  bounds on  the terms $|a_n(t)|$ for large values of $n$.    Our bounds will involve a number of variables in 
whose  choice we have  some limited flexibility.    In particular, we can select these so as to optimise the error terms.

\subsection{The contraction ratio $\theta$ for  expanding maps}
Let us assume that we can associate a Markov partition $\mathcal P = \{ P_1, \cdots, P_K\}$.
In the particular case that $T: X \to X$ is an expanding map we can consider 
\begin{enumerate}
\item
charts and the complexification of the maps $T$ 
to neighbourhoods $ \mathbb C^d \supset U_i \supset P_i$
($i=1, \cdots, K$);  and 
\item
 associated contractions $\psi_{ij}:U_i \to U_j$
wherever $T(\hbox{\rm int}(U_j)) \supset \hbox{\rm int}(U_i)$.
\end{enumerate} 

In more fortunate situations we can assume that we have \emph{Bernoulli} contractions where 
$U= U_i $ ($i=1, \cdots, K$) are identical and 
$\psi_{ij}$ = $\psi_i$ ($i=1, \cdots, K$) are independent of $j$.
(This applies in the case of Example \ref{doubling} and the Lanford map in Example \ref{lanford})

\begin{definition}
Choose $0 < \Theta_i < 1$ such that 
we can choose complex polydiscs
$ P_i \subset B(c_i, r_i) \subset U_i$
where 
$$
B(c_i, r_i) = \{\underline z \in \mathbb C
\hbox{ : } |z_i - c_i| < r_i, \quad i = 1, \cdots, K\}
$$
such that 
$$
\overline {\psi_{ij}(B(c_i, r_i))} \subset B(c_j, \Theta_j r_j).
$$
\bigskip
Let $\theta := \max_i \{\Theta_i^{1/K}\}$.  
In the particular case of Bernoulli contractions we can take 
$\theta := \max_i \{\Theta_i\}$.
\end{definition}

\medskip
\begin{rem}
It may not always  be possible to extend the contractions analytically to discs about elements of the 
partitions.  However, in that case we can consider instead the more  refined partition by elements $P_{i_0} \cap T^{-1} P_{i_1}  
\cap \cdots \cap T^{-n} P_{i_n} $, for suitable $n$.
\end{rem}
\medskip

\begin{example}
For an expanding map of the interval then we can take $d=1$.  
For a Bernoulli expanding map $T$ we would like to take $k=1$ (i.e.,  a single disc $B(c,r)$ and contractions 
$\psi_i: B(c,r) \to B(c,r)$ for $i=1, \cdots, k$) arising from the inverse branches of $T$.
The $P_i$ will be a partition of the interval into  subintervals and the discs 
$P_i \subset U_i \subset \mathbb C$ extend into the complex plane.
\end{example}

 \medskip
 \noindent
\emph{First parameter choice:} Choose  a real number $0 < \theta < 1$.  
 
 \bigskip
 There is no canonical choice  of $\theta$ and one can try to 
 arrange the partition and the polydiscs so as to minimise the possible choices.
 This can be achieved either by trial and error, or by simple calculus.

\begin{figure}
\centerline{
\begin{tikzpicture}
\draw (0,3.5) -- (0,-3.5);
\draw (3.5,0) -- (-3.5,0);
\draw[black] (0,0) circle (3cm);
\draw[black, dashed] (0,0) circle (2cm);
\draw[red] (1,0) circle (1cm);
\draw[red] (-0.5,0) circle (1.3cm);
\draw[<->] (0,0) -- (-2.12,2.12);
\node at (-1.7,1.9) {$r$};
\node at (0.2,0.2) {$x$};
\draw[<->] (0,0) -- (-1.42,-1.42);
\node at (-0.9,-0.5) {$\theta r$};
\node at (2.9,0.5) {$\psi_2(B(x,r))$};
\node at (-2.7,0.5) {$\psi_1(B(x,r))$};
\end{tikzpicture}
}
\caption{The choice of $0 < \theta < 1$ for two contractions $\psi_1, \psi_2: B(c,r) \to B(c,r) $ with $B(c, \theta r) \supset \psi_1(B(c,r)) \cup \psi_2(B(c,r))$.}
\end{figure}
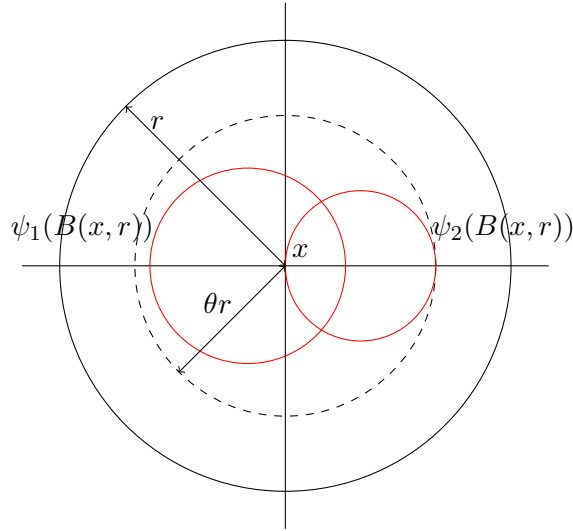

In the next subsection we begin to elaborate on the description of the error bounds mentioned 
in subsection \ref{illustrative}.

\subsection{Bounds on the determinant coefficients}
The key to obtaining validated estimates on characteristic values is to get accurate bounds
on the coefficients $a_n$, especially  for large $n \geq 1$.   

\bigskip
\noindent
 \emph{Second parameter choice:}
Choose a natural number  $N > 0$.

\bigskip
\noindent
This should  be chosen as large as is practicable.
Typically this will depend on the time, computer memory  and computing power available to compute periodic points.
For $\ell$ contractions, to compute $a_n$ for $1 \leq n \leq N$ would require making estimates on up to  of the order of 
$\ell^N$ periodic points.    This involves quantities $\beta_k$, $t_m$ and $B_k$ which we define below.

\begin{enumerate}
\item[i)]
Let   $\epsilon_1 = \epsilon_1(N) > 0$ be the corresponding error bound, 
by which we mean  that 
for $1 \le n\le N$
the coefficient  $a_n$
can be computed
  to a guaranteed  accuracy of no more than
 $\epsilon_1$.
 \end{enumerate}
 
 \bigskip
\noindent
 \emph{Third parameter choice:}
 Choose a natural number $L>N$.
 
 \bigskip
 
\noindent 
Typically this choice will depend on the time, computer memory  and computing power available to compute integrals.
This will involve us in numerically integrating (to appropriate precision) approximately $L-N$ integrals.

 \begin{enumerate}
 \item[ii)]
Let  
 $\epsilon_2 = \epsilon_2(L) > 0$ be the 
  bound on the tail 
$$\sum_{k=L}^\infty 
\|\mathcal L q_k\|^2  \leq \frac{C\theta^L}{1-\theta}< \epsilon_2$$
where 
$$
\begin{aligned}
\|\mathcal L q_k\|^2
&= r_i^{-2k}\int_{0}^{1}
| \sum_{j} \psi_j'(c_j + r_i e^{2\pi i t})(c_i + r_i \psi_j (e^{2\pi i t}))^k
|^2dt
\leq C\theta^n
\end{aligned}
$$
and where:
\begin{enumerate}
\item
$C = K^2\max_j \|\psi_j'\|_\infty^s$;
and 
\item
$q_k(z) = (z-c_i)^kr_i^{-k}$, with $i = i(k)$ (mod $K$),
\end{enumerate}
 where $K$ is the number of discs needed.
(In the Bernoulli case we have the integral around the same curve).
\end{enumerate}

\begin{enumerate}
\item[iii)]
 Let  $\epsilon_3 = \epsilon_3(L,N) > 0$ be the bound on  rigorous computational estimates $(\beta_k)_{k=N}^L$ such that 
$$
 \beta_k
 := \|\mathcal L q_k\|^2
= r_i^{-2k}\int_{0}^{1}
\left| \sum_{j} \psi_j'(c_j + r_i e^{2\pi i t})(c_i + r_i \psi_j (e^{2\pi i t}))^k
\right|^2dt
\hbox{ for $N <  k < L$}
$$
up to an error of $\epsilon_3$.

\end{enumerate}

We can now use these choices to define a sequence $(t_m)$ for the next set of bounds.

\begin{definition}
We can define a sequence of positive real numbers
$$
t_m:= 
\begin{cases}
\left(\sum_{k=m+1}^L
\beta_k + (L-m)\epsilon_3 + \epsilon_2 \right)^{1/2}
&\hbox{ for }
 m \leq L\\
 C \theta^m  &\hbox{ for } m > L.
 \end{cases}
 $$
for  $C>0$ as above.
\end{definition}

In particular, these numbers will tend to zero as $m\to\infty$.
We can combine  the values of $t_n$ by introducing the following definition of $B_k$.

\begin{definition}
For $1 \leq k\leq L$,  define  positive real numbers
$$
B_k :=  \sum_{m_1 < \cdots < m_k \leq L} 
t_{m_1} \cdots t_{m_k}.
$$
\end{definition}

\noindent
Typically, the coefficients $B_k$ will tend to zero quite quickly. 
Moreover, $B_k$ is defined  up to  an error 
$$
\epsilon_4 :=
\epsilon_3
\left(\max\{t_i\}\right)^{L-1}L^k.$$

The quantities $\beta_k$, $t_m$ and $B_k$ are used in giving 
the following bounds on the coefficients $a_n$ for the determinant.

\begin{thm}[Coefficient bounds]\label{coef}
Let $c = \frac{1}{\prod_{k=1}^\infty (1-\theta^j)}$.  Then
\begin{enumerate}
\item
for  $N < n \leq L $ we can bound,
$$
|a_n| \leq 
\gamma:=
c\sum_{k=1}^{n} (B_k +\epsilon_4)   (\theta^L C)^{n-k}; \hbox{ and }
$$
\item
for $n > L$   we can bound
$$
|a_n| 
\leq 
\xi  (\theta^L C)^{n}
\hbox{ where }
\xi:=
c
 \left(
\sum_{k=1}^{L} (B_k +\epsilon_4)  (\theta^L C)^{-k}\right).
$$
\end{enumerate}
\end{thm}

\noindent
In particular,  we get effective bounds for  $|a_n|$ when $n > N$
(combining these different bounds for $N < n \leq L$ and $n > L$).

\begin{rem}
In part 1 of Proposition  \ref{coef} we have used the simply proven bound
$$\sum_{r_1, \cdots, r_{n-k}=L}^\infty 
C^{n-k} \theta^{r_1+ \cdots +  r_{n-k}}
\leq 
c (\theta^L C)^{n-k}.
$$
\end{rem}

The proof of Theorem \ref{coef} follows the same lines as the arguments in 
\cite{jpeffective,jpzaremba,jpv}.  
For the convenience of the reader 
we give a brief account of 
the underlying operator theory ideas in the proof in Appendix A.

\subsection{Applying the bounds}
There are two different ways that the  bounds in Theorem \ref{coef} might be applied to estimate the accuracy of our approximations to the relevant characteristic values, depending on the quantity in question.

\bigskip

\noindent
{\bf (a) Explicit values.}
Assume that we have an expression that can be written in terms of $d(z,t)$ and its derivatives 
$\frac{\partial^{i+j}}{\partial z^i \partial s^j} d(z,s)$ ($i,j \geq 0$).  
The value of $d(z,t)$ can be approximated using the preceding estimates
$$
\left|
d(z,t)- \left(1+ \sum_{n=1}^N a_n(t) z^n  \right)
\right|
\leq  N |z|^N \epsilon_1 + \sum_{n=N+1}^M \gamma_n |z|^n + 
\xi     \frac{|z|(\theta^L C)^L}{1- |z|(\theta^L C)}.
$$
More generally,  we can bound 
$$
\left| \frac{\partial^{i+j}}{\partial z^i \partial s^j} d(z,s)
-
\sum_{n=1}^N z^n \frac{n!}{(n-i)!}
\frac{\partial^i}{\partial s^i} a_n(s)\right|
$$ 
where we bound $\frac{\partial^i}{\partial s^i} d(z,s)$ using Cauchy's theorem.

\bigskip
These estimates can be applied, for example, to computing the error terms in the estimates on 
Lyapunov exponents, variance and linear response as in subsection \ref{quantities}.
\bigskip

\noindent 
{\bf (b) Implicit values.}
If we are seeking a zero of $d(z,s)$
then provided we can choose $z_1 <  z_2$ close
and with validated bounds 
$$d(z_1, s) < 0 < d(z_2, s)
\hbox{  
or }
d(z_2, s) < 0 < d(z_1, s)$$
using $z=z_1$  or $z=z_2$,
   then there must be a zero in the interval $[z_0,z_1]$.

\bigskip
\noindent
This approach can be  used to estimate errors in computing  the rates of mixing and Hausdorff dimension as described in section \ref{mixing}.  We will illustrate this in a specific instance in the next section.

\section{A worked example:  Hausdorff dimension of the set $E_2$}
To illustrate how Theorem \ref{coef} can be applied in practice, we want to consider a particular concrete problem.
In particular, we will describe its use in estimating the Hausdorff dimension of a specific Cantor set (see \cite{jpeffective}).  

Recall from Example \ref{deleteddigitsexample} that $E_2$ is the subset of
$[0,1]$ consisting of those reals whose continued fraction expansion contains only the numbers 1 and 2.
In other words, if
$$
T_1(x) := \frac{1}{1+x}\ \hbox{ and }\  T_2(x) := \frac{1}{2+x}
$$
then $E_2$ is the corresponding
limit set (i.e.~the smallest non-empty closed set $X$ such that $T_1(X) \cup T_2(X) = X$).
Let us consider the estimation of the 
Hausdorff dimension of 
$E_2$, referring to \cite{jpeffective} for full details.

\bigskip
Defining $\underline i = (i_1, \cdots, i_n) \in \{1, 2\}^n$ and  $|\underline i|=n$,
and letting $x_{\underline i} = T_{\underline i}(x_{\underline i})$ be the fixed point for 
$$T_{\underline i} = 
T_{i_1} \circ \cdots \circ T_{i_n}: [0,1] \to [0,1]\,,$$
we have the   determinant 
$$
d(z,t):=
\exp \left(-
\sum_{n=1}^\infty \frac{z^n}{n}
\sum_{|\underline i|=n} 
\frac{|(T_{\underline i})'(x_{\underline i})|^t}
{1- (T_{\underline i})'(x_{\underline i})}
\right)\,,
$$  
and $\dim_H(E_2)$ is the value $t$ such that
 $d(1,t)=0$.
This value is then approximated
using the following four steps:

\begin{enumerate}
\item
For each  $t$ approximate  $ z \mapsto d (z, t)$ by a polynomial
 $z \mapsto d_N (z, t)$; 
\item Set $z=1$ and consider $t \mapsto d_N (z, 1)$;
\item Solve for $t_N = t$:  $d_N (1, t)=0$; 
\item
Then $t_N \to \hbox{\rm dim}_H(E_2)$ as $N \to +\infty$.
\end{enumerate}

In particular, for each $N \in \mathbb N$ we obtain an approximation 
$t_N$ to $\hbox{\rm dim}_H(E_2)$.
The sequence $t_N$ gives an intuitive estimate on the quality of the approximation in terms of the difference $|t_N - t_{N-1}|$.

We can write the series expansion
$$
d(z, t) = 1+ \sum_{n=1}^\infty a_n(t) z^n
= \underbrace{1+ \sum_{n=1}^N a_n(t) z^n}_{=: d_N(z,t)}
+ \underbrace{\sum_{n=N+1}^\infty a_n(t) z^n}_{=:\epsilon_N(z,t)}
$$
for some $N \geq 1$, and
take for the approximating polynomial
$$
d_N(z,t)=
1+ \sum_{n=1}^N a_n(t) z^n
$$
where $N$ is chosen
to be sufficiently large that 
(with $z=1$ and $0 \leq t \leq 1$) the  error $\epsilon_N$ is small, 
but
sufficiently small that the  terms $a_n(t)$, $n=1,2, \cdots, N$ can be calculated in a reasonable time.
In the present setting one might  choose $N=25$, say.

\smallskip
We can explain part of the mechanism for effective estimates on the coefficients $a_n$ ($n \geq 1$)
as follows.

\bigskip
\noindent
\emph{Step 1.}  
Choose $z_0 \in \mathbb R$ and $r>0$ such that
$$D = \{ z \in \mathbb C \hbox{ : } |z-z_0| < r\} \supset [0,1]
\hbox{ and } 
T_1D, T_2 D \subset D.$$
For example, we could let $z_0 = 1$ and $r = \frac{3}{2}$.
Consider the family of  transfer operators defined 
 on  analytic functions $w: D \to \mathbb C$ by 
$$
{\mathcal L}_t w(z) = \frac{1}{(z+1)^{2t}} w\left(
\frac{1}{z+1}
\right)
+ \frac{1}{(z+2)^{2t}} w\left(
\frac{1}{z+2}
\right), \quad t \in \mathbb R.
$$

\bigskip
\noindent
\emph{Step 2.}   
Let 
$q_k(z) = \frac{(z-z_0)^k}{r^k}$, for   $z\in D$ and $k \geq 1$.
We can then define
$$
t_m = \left(\sum_{k=m-1}^\infty \|\mathcal L (q_k)\|^2 \right)^\frac{1}{2}
(m \geq 1)
\hbox{ where }
\|q_k\|^2 = \int_0^1| q_k(z_0 + re^{2\pi i t})|^2 dt  \quad (k \geq 1).
$$

\bigskip
\noindent
\emph{Step 3.} 
We can  bound the coefficients $a_n$  ($n > N$) by
$$
|a_n|  \leq \sum_{m_1 < \cdots < m_n} t_{m_1} t_{m_2} \cdots t_{m_n}.
$$

We will give more details of the underlying operator theory in Appendix A.

Given $M > 0$ (in the present setting one can choose $M=600$) we can numerically estimate $\|\mathcal L(q_k)\|$ for $k \leq M$, 
and we can trivially bound $\|\mathcal L(q_k)\|$ for $k > M$.
Combining these various bounds gives the estimate.

\section{Future directions}
In this note we have discussed how the periodic point algorithm works and  how the error terms can be efficiently estimated.
 However, there may be further scope to fine-tune the underlying analysis and estimates, 
  perhaps by using different Hilbert spaces or other operators.  This leads to the first  very general question.

\begin{question}
Can we improve the approach to computing the dynamical determinant and, in  particular,  the error estimates? 
\end{question}

We have illustrated the general  approach with a number of  applications and examples.  However, we now want to propose 
 some further potential applications.

  \subsection{Dynamical invariants and Lyapunov exponents}
  
  We have already discussed the theoretical use of our method to rigorously estimate certain dynamical quantities.
  However this has only been done practically in only a small number of cases (e.g., Lyapunov exponents, variance).  
  However, it remains a task yet to be  completed to explore how to estimate rigorously other quantities (e.g., resonances, linear response) for simple test examples (such as the Lanford map).
  
  \begin{question}
  Apply this approach to compute more dynamical invariants in simple examples.
  \end{question}

 Lyapunov exponents also occur naturally in the theory of 
 random matrix products, which has been  an  active area of research since the 
 pioneering work of 
Kesten and Furstenberg in the 1960s.  For example, given $k \times k$ square  matrices $A_1, \cdots, A_k$
{\it with positive entries}
we define the {\it Lyapunov exponent} by
$$
\lambda = \lim_{n\to +\infty} \frac{1}{d^n}
\sum_{i_1, \cdots, i_n \in \{1, \cdots, d\}}
\frac{\log \| A_{i_1} \cdots A_{i_n}\|}{n}.
$$
There was an  implementation of the basic algorithm in 
\cite{pollicottinventiones},

\begin{question}
Can the error bounds in the computation of $\lambda$ be made more effective?
\end{question}
 
 An interesting application is to the case of binary symmetric channels in information theory. 
  In this context there are positive matrices and $\lambda$ is related to a useful value called the 
   Entropy Rate.  
 	
\subsection{ Connections to number theory}\label{numbertheorysubsection}
We have already mentioned the application to the density one Zaremba 
conjecture (see \cite{jpzaremba}).  However, we now want to describe a different application to number theory.

Given an irrational number $\alpha$ we can associate the number
$$
\mu(\alpha) = \limsup_{p,q \to +\infty}|q^2| \left|\alpha - \frac{p}{q}\right|
$$
(i.e., the best constant in diophantine approximation for $\alpha$).  The {\it Lagrange spectrum} is defined to be the set
$\mathcal L = \{
1/\mu(\alpha) \hbox{ : }  
\alpha \in \mathbb R - \mathbb Q \}$.
On the other hand, we can consider those binary quadratic forms
$f(x,y) = a x^2 + b xy + c y^2$ ($a,b,c  \in \mathbb R$) with discriminant $D(f) = b^2 - 4ac > 0$  and denote 
$$\lambda(f) =  \inf \left\{|f(x,y)| \mbox{ : } (x,y \in \mathbb Z^2 -\{(0,0)\})\right\}/\sqrt{D(f)}.$$   The {\it Markov spectrum}  is defined to be the set 
$\mathcal M = \{
1/\lambda(f) \mbox{ : }  
\alpha \in \mathbb R - \mathbb Q \}$.
 It is known that $\mathcal L \subset \mathcal M$ 
and
 Matheus \& Moreira showed that 
 $0.513 \cdots  < \dim_H(\mathcal M - \mathcal L) < 0.98 \cdots$ \cite{matheusmoreira2}.  Moreover, in their article they conjecture  that the bounds can be  improved to $\dim_H(\mathcal M - \mathcal L) < 0.88$, based on empirical estimates using our algorithm.

\bigskip 
 \begin{question}
Obtain  improved rigorous bounds on  $dim_H(\mathcal M - \mathcal L)$.
 \end{question}

\medskip 

This involves rigorously computing the Hausdorff dimension of limit sets associated to iterated function systems 
$\{\phi_i: I \to I\}$,
but with a Markov condition, i.e., there is a $0$-$1$ matrix $A$ and compositions $\phi_i \circ \phi_j$ are only allowed if $A(i,j)=1$.
Whereas the basic algorithm still applies in this setting, the major complication is to get effective error estimates.

\medskip

\subsection{Spectral geometry and the Selberg zeta function}\label{geometrysubsection}

Given a  compact surface  $V$, with a metric $\rho$ of  constant negative curvature, we can associate the 
 Selberg zeta function  defined by 
$$
Z_\rho(s) =
 \prod_{n=0}^\infty 
 \prod_\gamma\left(1 - e^{-s(n+s)l(\gamma)} \right), \quad s \in \mathbb C,
$$
where $\gamma$ denotes a closed geodesic of length $l(\gamma)$.
To relate this to our analysis, we recall that we can associate a piecewise $C^\omega$ expanding map of the circle 
(using the Bowen-Series approach) and then the  zeta function can be written in terms of  the determinants 
$\det(I-\mathcal L_s)$
of the associated transfer operators $\mathcal L_s$
(see e.g.~\cite{mayercmp, ruelleinventiones}).

The zeros of $Z_\rho(s)$ have a spectral interpretation, in terms of eigenvalues of the Laplacian on $(V, \rho)$, 
but these can be well estimated using other techniques.
Other special values such as 
$Z_\rho'(0)$ can be written in terms $\frac{d}{ds}\det(I-\mathcal L_s)|_{s=0}$ 
and this is  proportional to the much studied
{\it  determinant of the Laplacian},  originally defined in terms of the spectrum of the Laplacian
(see e.g.~\cite{dhokerphong, friedinventiones, sarnakcmp}).

\begin{question}
Can we get useful  bounds on $Z_\rho'(0)$ ?
\end{question}

There is a well-known problem of Sarnak to show that there is a (local) minimum for the determinant that occurs at very symmetric hypergeometric surfaces which could be addressed using this approach.

In a different direction, the  Weil-Petersson metric 
is a classical distance  on 
the space of  such Riemann metrics $\rho$.  There is a particularly useful  thermodynamic interpretation for the Weil-Petersson metric 
due to McMullen in
 terms of the second derivative of the associated thermodynamic pressure in the context  above of the piecewise $C^\omega$ expanding maps.  
 
 \begin{question} 
 Can one get effective estimates on the  Weil-Petersson metric?
 \end{question}
 
  This could then be used to explore empirically  Weil-Petersson metric on the space of metrics.
 Moreover, there are higher dimensional analogues of Weil-Petersson metrics  in
 \cite{bcls} which could be similarly analysed.

\appendix
\section*{Appendices}
\addcontentsline{toc}{section}{Appendices}
\renewcommand{\thesubsection}{\Alph{subsection}}


\subsection{Appendix: Operator theory and coefficient bounds}
We now give a little more detail on how the bounds on the coefficients $a_n$ arise  using  operator theory.

\subsubsection{Operator theory}
Assume for notational convenience that $T:X \to X$ is a Bernoulli expanding map of the unit interval $X$.
Assume that we can choose an open disc $X \subset U \subset \mathbb C$ so that 
$T$ extends analytically to $U$. 
In the particular case that $D= \{z \in \mathbb C \hbox{ : } |z - z_0| < r\}$ then 
we can let $\mathcal H$ be the Hardy Hilbert space 
of analytic functions on $U$ for which the norm $\| \cdot \|$ is given by
$$\|f\|^2 = \sup_{\rho < r} \int_0^1| f(z_0 + \rho e^{2\pi i t})|^2 dt <\infty, \quad f\in \mathcal H.$$

\medskip
\noindent 
{\bf Example A.1.1} 
If   $D = \left\{z \in \mathbb C \hbox{ : } |z - \frac{1}{2}| < 1\right\}$ then $f \in \mathcal H$  presented in the form   
$$f(z) = \sum_{n=0}^\infty b_n \left(z - \frac{1}{2}\right)^n$$
 has norm
$\|f\|^2 = \sum_{n=0}^\infty |b_n|^2$.
\medskip

Let $\mathcal L = \mathcal L_{G,T}: \mathcal H \to \mathcal H$ be a transfer operator.
We can then assume that $\mathcal L$
is a trace class operator and a simple Lidskii-type identity
relates periodic orbits to the spectrum
$$
\hbox{\rm trace} (\mathcal L^n) = 
\sum_{T^nx=x} 
\frac{
	\exp\left(\sum_{i=0}^{n-1} G(T^ix)\right)
}{1 - 1/(T^n)'(x)}
$$
and then 
$$D(z)
= \det(I - z\mathcal L) = \exp \left(
-\sum_{n=1}^\infty \frac{z^n}{n} \hbox{\rm trace} (\mathcal L^n)\right)$$
where both sides depend on $G$ and $T$ \cite{robert}.
The basic idea is  to use the operator description of $D(z)$ to get estimates on the $a_n$ via the
approximation numbers $s_l$ ($l \geq 1$) defined by 
$$s_l := \inf \left\{\|\mathcal L - K\| \hbox{ : } K \hbox{ bounded linear operator} \hbox{ of finite rank }   l\right\}, \quad l \geq 1\,.$$
For definiteness, let $(e_n)_{n=0}^\infty$ be a complete orthonormal family for $\mathcal H$.
	For any
	$$f = \sum_{n=0}^\infty c_n e_n
	\in \mathcal H  \hbox{  with } \|f\| = \sqrt{\sum_{n=0}^\infty |c_n|^2}=1$$
	 and  each $l \geq 1$ we can define
	$K_lf = \sum_{n=0}^{l-1} c_n \mathcal L(e_n) \in \mathcal H$.
In particular, we can bound 
$$
\begin{aligned}
s_l
\leq &\|(\mathcal L - K_l)f\|
\leq \sum_{n=l}^\infty |c_n|\|\mathcal L(e_n)\|
\cr
&\leq \sqrt{\sum_{l=n}^\infty |c_n|^2}
\sqrt{\sum_{n=l}^\infty \|\mathcal L(e_n)\|^2}
\leq 
\sqrt{\sum_{n=l}^\infty \|\mathcal L(e_n)\|^2}
\end{aligned} 
\eqno(A1)
$$
using the Cauchy-Schwarz inequality.

\medskip
\noindent 
{\bf Example A.1.2} (Example A.1.1 revisited).
If  $D = \{z \in \mathbb C \hbox{ : } |z - \frac{1}{2}| < 1\}$
then we can take  $e_n(z) = (z-\frac{1}{2})^n/\sqrt{n}$, for $n \geq 0$.
\medskip

We can now  bound the coefficients $a_n$ ($n > N$) by
$$
|a_n|  \leq \sum_{m_1 < \cdots < m_n} s_{m_1} s_{m_2} \cdots s_{m_n}, \eqno(A2)
$$
where the summation is over distinct natural numbers $m_1, \cdots, m_n \in \mathbb N$
such that $m_1 < \cdots < m_n$, 
using \cite[Cor.~VI.2.6]{ggk}.

Given $M > 0$ (in the present setting one can choose $M=600$) we can numerically estimate $\|\mathcal L(e_k)\|$ for $k \leq M$, 
and we can trivially bound $\|\mathcal L(e_k)\|$ for $k > M$.



\subsubsection{Basic convergence}
The basic convergence (which follows from  work of Ruelle \cite{ruelleinventiones}, after Grothendieck) only requires the following bound (referred to in \cite{jpeffective,jpzaremba,jpv} as an {\it Euler bound}, after \cite{euler}):

\bigskip
\noindent
		If there exists 
		 $C >0$ and $0<r<1$ with
		$s_l \leq Cr^l$  for $l \geq 1$ then 
		$$|a_n| \leq 
		\frac{C^n r^{n(n+1)/2}}{(1-r)(1-r^2) \cdots (1-r^n)} = O(r^{n(n+1)/2})\eqno(A3)$$
		for $n \geq 1$.
\bigskip	

	In particular, the coefficients $a_n$ tend to zero fast enough to
	ensure absolute convergence of  the series
	describing the various quantities of interest
	 (cf. Part 1 of Corollary \ref{cor} and Corollary \ref{super}).

\begin{rem}
We can now see the technical limitations of this method.  
In (A3) we have a bound $|a_n| = O(r^{n^2})$.  However, if we look at $d$-dimensional analytic maps 
then we have a slightly weaker  bound $|a_n| = O(r^{n^{1+\frac{1}{d}}})$ which will make the approach marginally less efficient.    Moreover, if there are infinitely many branches then typically we lose the advantages of this method because the terms $a_n$ involve infinite summations which need to be approximated individually.
\end{rem}

	\subsubsection{Better estimates}
	Assume that we can compute numerically $a_1, \cdots, a_N$ for some $N$.
	We could  then  use the above bound (A3)  to bound the terms $a_n$ for $n > N$
	 (although this wasn't the original purpose of this approach nor is it necessarily a very efficient  approach).   
	 The basis of the improved approach involves:
	\begin{enumerate}
		\item 
		again computing  $a_1, a_2, \cdots, a_N$ using the periodic points of period at most $N$ (where $N$ is chosen depending on the limitations of our computer).
		\item  bounding $|a_n|$, for $n\geq N$ using (A2) and different bounds on the $s_l$: 
		\begin{enumerate}
			\item
			for $N < l \leq  M$ we can 
			bound 
			$s_l$ using (A1) (where $M$ is chosen depending on the limitations of our computer); and
				\item 
			For $M < l$
			we can use the  bound 
			$s_l \leq C r^l$.
		\end{enumerate}
	\end{enumerate}
	
	There are further refinements possible, but this explains the gist of the approach (see \cite{jpeffective}).
	
	\subsubsection{Generalisations}
	Finally, we can describe a plan for how this method can  be extended to Anosov diffeomorphisms and flows.
	
	If we have a $C^\omega$ Anosov diffeomorphism $T$  then it is more appropriate to write
	$$
	D(z)  =  D_{G,T}(z) = \exp\left(- \sum_{n=1}^\infty 
	\frac{z^n}{n} \sum_{T^nx=x} 
	\frac{\exp\left({\sum_{i=0}^{n-1}G(T^ix)}\right)}{\det(I - (DT)^{-1}(x))}
	\right), \quad z\in \mathbb C,
	$$
	An extra feature here is that we 
	only need to take $G_t(x) = t g(x)$ because of the contributions from $\det(I - (DT)^{-1}(x))$.
	
	One can adapt this approach using either the less fashionable device of Markov partitions and the approach of Rugh, or using anisotropic spaces of functions when applicable (for example on the torus in the context of Faure-Roy).  The key point is simply to find a setting to which the functional analysis applies.
	
	For $C^\omega$ Anosov flows it is perhaps  simpler to use Markov Poincar\'e 
	sections to accommodate the operator theory being used.


\subsection[Appendix: Asymptotic expansion for quadratic Julia sets]{Appendix: Asymptotic expansion for the dimension of quadratic Julia sets}
	
As mentioned in Section \ref{repellerssubsection}, we detail here	
a derivation of the asymptotic expansion
$\dim_H(J_c) = 1 + \frac{|c|^2}{2\log 2} + O(|c|^3)$ s
from Theorem \ref{ruelleasymp}.
We will follow the approach of Zinsmeister in \cite{Zin00}, taking the opportunity to correct a few minor errors.
\footnote{We are  grateful to Michel Zinsmeister for useful discussions on these points.}
We can write  $c = x+iy$ and define 
\begin{equation*}
F(d,x,y) := P\left(-d(c)\log(2|\Psi_{x+iy}(z)|)|_{J_0}\right).
\end{equation*}
We can then solve  for $d(x,y) = \dim_H(J_c)$ satisfying $F(d(c),x,y)=0$.
Expanding $d(c)$ as a Taylor series around $c=0$, we obtain:
\begin{equation*}
d(x,y) = 1 + \alpha x + \beta y + \gamma x^2 + \delta y^2 + \epsilon xy + o(x^2+y^2).
\end{equation*}
For $c$ close to $0$ we have that $\dim_H(J_c) \geq \dim(J_0) =1$ and 
it  follows that $\alpha$ and $\beta$ must be zero. 
By expanding $F$ as a series around the point $(1,0,0)$, we get:
\begin{align*}
0 = F(d,x,y) = x F_x + yF_y &+x^2\left(\delta F_d + \frac{1}{2}F_{xx}\right) +y^2\left(\gamma F_d + \frac{1}{2}F_{yy}\right)\\
&+ xy \left(\epsilon F_d + F_{xy}\right)
+ o(x^2+y^2)
\end{align*}
where $F_x $, $F_y$, $F_{xx}$, $F_{yy}$ and $F_{xy}$ are the partial derivatives as $x=d$ and $y=0$, 
and thus 
\begin{equation*}
F_x = F_y = 0, \, \, \gamma = -\frac{F_{xx}}{2F_d}, \, \, \delta =-\frac{F_{yy}}{2F_d}, \, \, \epsilon = -\frac{F_{xy}}{F_d}.
\end{equation*}

To begin, since $F(d,0,0) = P(-d\log2)  =  \log2 - d \log 2$ we have $F_d = -\log 2$.
Next, we can  expand  $\Psi_c$ in terms of  $c$:
$$
	\Psi_c(z) = z + a_1(z)c +a_2(z)c^2 + o(c^2). \eqno(B1) 
$$
 Since $\Psi_c$ satisfies $\Psi_c \circ f_0 = f_c \circ \Psi_c$, we can write 
$\Psi_c(z^2) = (\Psi_c(z))^2 + c$.  Substituting the expansion (B1)  and  comparing the coefficients for $c$ and $c^2$, we obtain\footnote{This is the first of the typographical inaccuracies in \cite{Zin00}, replacing $(a_1(z))^2$ with $a_1(z^2)$.} 
$$
a_1(z^2) = 2za_1(z) + 1 \, \, \, \mathrm{and} \, \, \, (a_1(z))^2 + 2za_2(z) = a_2(z^2).
\eqno(B2) 
$$
The first equation in  (B2) has a series solution:
\footnote{This is the second of the minor changes to  \cite{Zin00}, replacing 
$z^{-2^n}$ with $z^{-2^{n+1}}$}
$$
a_1(z) = -z\sum_{n=0}^{\infty}\frac{1}{2^{n+1}} z^{-2^{n+1}}.
$$
We next  write $a_2(z) = \sum_{n=-\infty}^{\infty} d_nz^n$ and compare coefficients in the second equation in (B2). In particular, we  obtain $d_n = 0$ for all $n \geq 0$, and $d_{-1}=d_{-2}=-1/8$. For $|z|= 1$ we have:
\begin{equation*}
|\Psi_c(z)| = |1 + b(z)c + \overline{z}a_2(z)c^2 + \ldots|,
\end{equation*}
where we write $b(z):=\overline{z}a_1(z)$ to make notation easier. 
Writing  $|\Psi_c(z)|^2 = (\mathrm{Re}(\Psi_c(z)))^2 +\mathrm{Im}(\Psi_c(z))^2$, we  have that 
\begin{align*}
\log|\Psi_c(z)| &= \log \left((1+\mathrm{Re}(b(z)c) + \mathrm{Re}(\overline{z}a_2(z))^2 + (\mathrm{Im}(b(z)c)))^2 + o(|c|^2)\right)^{1/2} \\
& = \frac{1}{2}\log \left( 1 +  2\mathrm{Re}(b(z)c) + |b(z)c|^2 + 2\mathrm{Re}(\overline{z}a_2(z)c^2) + o(|c|^2) \right).
\end{align*}
Using the   series $\log(1+w) = w - w^2/2 +w^3/3 + \ldots$ yields:
\begin{equation*}
\log |\Psi_c(z)| = \mathrm{Re}(b(z)c) +\frac{1}{2} |b(z)c|^2 + \mathrm{Re}(\overline{z}a_2c^2) - (\mathrm{Re}(b(z)c))^2 + o(|c|^2)
\end{equation*}
and then 
writing  $d=1+h$:
\begin{align*}
F(1+h,x,y) &= P(-(1+h)\log 2 -(1+h)\log|\Psi_c|) \\
&=P \left( -\log2 - \mathrm{Re}(b(z)c) - \frac{1}{2} |b(z)c|^2 - \mathrm{Re}(\overline{z}a_2c^2) + \mathrm{Re}(b(z)c^2) - h\log2 -\mathrm{Re}(b(z)c) + \ldots \right)\\
&= 0 + P_{g}'(-\log 2) + \frac{1}{2} P_{g,g}''(-\log 2) + o(|c|^2 + h^2),
\end{align*}
where 
$
g = -\mathrm{Re}(b(z)c) - \frac{1}{2} |b(z)c|^2 - \mathrm{Re}(\overline{z}a_2c^2) + (\mathrm{Re}(b(z)c))^2 - h\log 2 - h \mathrm{Re}(b(z)c).
$
We now compute $\int_{J_0} g \, d\lambda$, where $\lambda$ is the normalized Haar measure on $J_0$. We begin by calculating  the first term
\begin{align*}
	\int_{J_0}\mathrm{Re}(bc) \, d\lambda &
	= \frac{x}{2\pi} \int_{0}^{2\pi}  \sum_{n=0}^{\infty} \frac{\cos(2^{n+1}2\pi t)}{2^{n+1}} \, dt 
	- \frac{y}{2\pi} \int_{0}^{2\pi}  \sum_{n=0}^{\infty} \frac{\sin(2^{n+1}2\pi t)}{2^{n+1}} \, dt 
	= 0.
\end{align*}
Similarly
we can show that:
\begin{equation*}
	\int_{J_0} \mathrm{Re}(\overline{z}a_2(z)c) \, d\lambda = 0.
\end{equation*}
We briefly describe how to compute the remaining  values. We have\footnote{This is another slight  deviation from \cite{Zin00}, where the integral is computed to be $2(x^2+y^2)/3$ rather than $(x^2+y^2)/6$.}: 
\begin{align*}
\int_{J_0} (\mathrm{Re}(b(z)c ))^2 \, d\lambda 
&= \int_{J_0} \mathrm{Re}\left(
(x+iy)
\sum_{n=0}^{\infty} \frac{1}{2^{n+1}} 
\left(z^{-2^{n+1}}\right)\right)^2 \, d\lambda \cr
& =
\frac{x^2}{2\pi}  \int_{0}^{2\pi}
\left( 
 \sum_{n=0}^{\infty} 
 \frac{\cos(2^{n+1}2\pi t)}{2^{2(n+1)}}
 \right)^2
  \, dt 
+ 
\frac{y^2}{2\pi}  \int_{0}^{2\pi}
\left( 
 \sum_{n=0}^{\infty} 
 \frac{\sin(2^{n+1}2\pi t)}{2^{2(n+1)}}
 \right)^2
  \, dt 
\\
&= (x^2+y^2)\sum_{n=0}^{\infty} 
\frac{1}{2^{2n+3}}
= \frac{x^2 + y^2}{6}.
\end{align*}
 We also compute:
\begin{align*}
\int_{J_0} |b(z)|^2 \, d\lambda &= \int_{J_0} \left( \sum_{n=0}^{\infty} \left(\frac{1}{2^{n+1}} z^{-2^{n+1}}\right) \sum_{n=0}^{\infty} \left(\frac{1}{2^{n+1}} \overline{z}^{-2^{n+1}}\right) \right) \, d\lambda  = \int \sum_{n=1}^{\infty} \frac{1}{4^n} \, d\lambda = \frac{1}{3}.
\end{align*}
Therefore, $\frac{1}{2}\int |b(z) c|^2 \, d\lambda = |c|^2/6$ and we conclude:
$
\int_{J_0} g \, d\lambda = - h \log 2.
$
Writing $b(z) =b_1+ib_2$ gives:
\begin{equation*}
F(1+h,x,y) = -h \log 2 +x^2\left(\lim_{m\to \infty} \frac{1}{N} \int_{J_0} (b_1^N)^2 \, d \lambda \right) + y^2\left(\lim_{N\to \infty} \frac{1}{N} \int_{J_0} (b_2^N)^2 \, d \lambda \right) + \ldots,
\end{equation*}
where $b_i^N(z) = \sum_{k=0}^{N-1} b_i \circ f^k(z)$.  Calculations  similar to that for 
$\int (\mathrm{Re}(b(z) c))^2 d\lambda$ 
give  that:
\begin{align*}
\int_{J_0} (b_1^N(z))^2 \, d\lambda &
=
\frac{1}{2\pi}\int_{0}^{2\pi} \left(
 \sum_{k=0}^{N-1}
\sum_{n=0}^{\infty} \frac{1}{2^{n+1}} \cos\left(2 \pi t(2^{k}-2^{n+1})\right)\right)^2 dt
= \sum_{k=0}^{N-1} \frac{1}{4} = \frac{N}{4}
\end{align*}
and similarly $\int_{J_0} (b_2^N(z))^2 \, d\lambda = N/2$.
Therefore, since $F_{xx} = F_{yy} = \frac{1}{2}$ and $F_{xy} = 0$ we can write 
\begin{equation*}
F(1+h,x,y) = -h \log 2 + \frac{1}{4} x^2 + \frac{1}{4}y^2 + o(x^2 +y^2 +h^2).
\end{equation*}
So we obtain that $\gamma = \delta = \frac{1}{4\log 2}$ and $\epsilon = 0$. Therefore:
\begin{equation*}
d(c) = 1 + \frac{|c|^2}{4 \log 2} + o(|c|^2),
\end{equation*}
as claimed.


\begin{thebibliography}{111}



\bibitem{BZ}
O.~Bodart \& M.~Zinsmeister,
Dimension de Hausdorff des ensembles de Julia.
{\it Fund.\ Math.} {\bf 151} (1996) 121--137


\bibitem{bourgainkontorovich}
J. Bourgain \& A. Kontorovich,  
On Zaremba's conjecture,
{\it Ann. of Math.}, {\bf 180} (2014), 137--196.



\bibitem{bowenbook}
 R.~Bowen:  
 {\it Equilibrium states and the ergodic theory of Anosov
 diffeomorphisms.} 
 Springer Lecture Notes in Mathematics
 {\bf 470}
 Berlin, Springer (1975)

\bibitem{bowenpublihes}
R.~Bowen:
 Hausdorff dimension of quasi-circles.
{\it Publ.\ Math.\ (IHES)}
 {\bf 50} (1979) 11--25

\bibitem{bcls}
M. Bridgeman, 
R. D. Canary, 
F. Labourie \&
A. Sambarino, 
The pressure metric for Anosov representations,
{\it Geom. Funct. Anal.}, {\bf 25} (2015), 1089--1179.

\bibitem{bumby1}
R. T. Bumby,
Hausdorff dimensions of Cantor sets,
{\it J. Reine Angew. Math.},
{\bf 331}
 (1982), 192--206.


\bibitem{cvitanovic}
P. Cvitanovic, R. Artuso, R. Mainieri, G. Tanner and G. Vattay,
{\it Chaos: Quantum and Classical}
chaosbook.org/version12/chapters/ChaosBook.pdf


\bibitem{dhokerphong}
E.D'Hoker \& D. Phong,
 On determinants of Laplacians on Riemann surfaces,
 {\it Comm. Math. Phys.}, {\bf 104} (1986), 537--545.

\bibitem{dsz}
A. Douady, P. Sentenac \& M. Zinsmeister, 
Implosion parabolique et dimension de Hausdorff,
{\it C. R. Acad. Sci. Paris Sér. I Math.}, {\bf 325} (1997), 765--772.

\bibitem{euler}
L. Euler, {\it Introductio in Analysin Infinitorum}, Marcum-Michaelem Bousquet, Lausannae, 1748.

\bibitem{falconer}
K. Falconer
{\it Fractal geometry. Mathematical foundations and
applications}, 
John Wiley \& Sons, Ltd., 1990.

\bibitem{falknussbaum} R. S. Falk \& R. D. Nussbaum,
$C^m$ eigenfunctions of Perron-Frobenius operators and a new approach to numerical computation of Hausdorff dimension: applications in $\mathbb{R}^1$,
{\it J. Frac. Geom.}, {\bf 5} (2018), 279--337.

\bibitem{falknussbaum1} 
R. S. Falk \& R. D. Nussbaum, 
A New Approach to Numerical Computation of Hausdorff Dimension of Iterated Function Systems: Applications to Complex Continued Fractions,
Integral Equations and Operator Theory, 90 (2018)


\bibitem{fredholm} I. Fredholm, 
Sur une classe d?\'equations fonctionnelles,
{\it Acta Math.}, {\bf 27} (1903), 365--390.








\bibitem{fried} 
 D.~Fried: 
Zeta functions of Ruelle and Selberg I. 
{\it Ann.\ Sci.\ Ec.\ Norm.\ Sup.} {\bf 9} (1986) 491--517

\bibitem{friedinventiones}
D. Fried,
Analytic torsion and closed geodesics on hyperbolic manifolds,
{\it Invent. Math.},
{\bf 84} (1986),
523--540.




\bibitem{ggk}
I. Gohberg, S. Goldberg \& M. A. Kaashoek,
{\it Classes of linear operators vol.~1}, 1990,
Birkh\"auser, Berlin.

\bibitem{good}
I. J. Good,
The fractional dimensional theory of continued fractions,
{\it Proc. Camb. Phil. Soc.},
{\bf 37} (1941),
199--228.

\bibitem{grothendieckthesis}
A. Grothendieck,
Produits tensoriels topologiques et espaces nucl\'eaires,
{\it Mem. Amer. Math. Soc.},
{\bf 16}, 1955.

\bibitem{grothendieckbullsmf}
A. Grothendieck,
La th\'eorie de Fredholm,
{\it Bull. Soc. Math. France},
{\bf 84} (1956),
319--384.


\bibitem{holliday}
T. Holliday, P. Glynn and A. Goldsmith,
Shannon maats Lyapunov: Connections between information theory 
and dynamical  systems, Proc. 44th IEEE conference on Decision and Control, and 
European Control Conference 2005. pp. 1756-1763, Elsevier, 2005

\bibitem{hej}
D.  A. Hejhal, Eigenvalues of the Laplacian for$PSL(2, \mathbb Z)$: some new results and computational techniques, International Symposium in Memory of Hua Loo Keng, Vol. I (Beijing, 1988) Springer, Berlin, 1991, 59-102

\bibitem{hensley1989}
D. Hensley,
The Hausdorff dimensions of some continued fraction Cantor sets,
{\it J. Number Theory}, {\bf 33} (1989), 182--198.

\bibitem{hensley1996}
D. Hensley, 
A polynomial time algorithm for the Hausdorff dimension of continued fraction Cantor sets,
{\it J. Number Theory}, {\bf 58} (1996), 9--45.


\bibitem{huang}
S. Huang, An improvement to Zaremba's conjecture, 
{\it Geom. Funct. Anal.},
{\bf 25} (2015), 860--914.

\bibitem{jarnik}
I. Jarnik, Zur metrischen Theorie der diophantischen Approximationen,
{\it Prace Mat.-Fiz.}, {\bf 36} (1928), 91--106.

\bibitem{jpeffective} 
O. Jenkinson \& M. Pollicott,  
Rigorous effective bounds on the Hausdorff dimension of continued fraction Cantor sets: a hundred decimal digits for the dimension of $E_2$,
{\it Adv. Math.}, {\bf 325} (2018), 87--115.

\bibitem{jpzaremba}
O. Jenkinson \& M. Pollicott,  
Rigorous dimension estimates for 
 Cantor sets arising in Zaremba theory,
 {\it Contemp. Math.}, to appear.

\bibitem{jpv} O. Jenkinson, M. Pollicott \& P. Vytnova,
Rigorous computation of diffusion coefficients for expanding maps,
{\it J. Stat. Phys.},
{\bf 170} (2018), 221--253.


\bibitem{kh}
A. Katok and B. Hasselblatt,
Introduction to the modern theory of dynamical systems, C.U.P., Cambridge, 1995.
 

\bibitem{kontorovich} A. Kontorovich,
From Apollonius to Zaremba:
local-global phenomena in thin orbits,
{\it Bull. Amer. Math. Soc.},
{\bf 50} (2013), 187--228.


\bibitem{lanford} O. E. Lanford III, Informal remarks on the orbit structure of discrete approximations to chaotic maps,
{\it Exp. Math.}, {\bf 7} (1998), 317--324.

\bibitem{liverani}
C.~Liverani:
Rigorous numerical investigation of the statistical properties of piecewise expanding maps. A feasibility study,
{\it Nonlinearity} {\bf 14} (2001), 463--490. 




\bibitem{mcmullen3}
C. McMullen
Hausdorff dimension and conformal dynamics. III.
 Computation of dimension
{\it Amer. J. Math.}, 
{\bf 120} (1998), 691--721.
 

\bibitem{markov} A. Markov, Sur les formes  binaires ind\'efinies,
{\it Math. Ann.}, {\bf 15} (1879), 381--406.

\bibitem{matheusmoreira} C. Matheus \& C. G. Moreira,
$HD(M\setminus L) > 0.353$,
{\it arXiv preprint, arXiv:1703.04302}.

\bibitem{matheusmoreira2} C. Matheus \& C. G. Moreira,
$HD(M\setminus L) < 0.986927$,
{\it arXiv preprint, arXiv:1708.06258}.


\bibitem{mu}
R. Daniel Mauldin and Mariusz Urba\'nski, Dimensions and measures in infinite iterated
function systems, Proc. London Math. Soc. (3) 73 (1996), no. 1, 105-154. 

\bibitem{mayercmp} D.~H.~Mayer,
On the Thermodynamic Formalism for the Gauss Map,
{\it Comm. Math. Phys.}, {\bf 130} (1990), 311--333. 







\bibitem{pollicottinventiones}
M. Pollicott, Maximal Lyapunov exponents for random matrix products,
{\it Invent. Math.}, {\bf 181} (2010), 209--226.

\bibitem{pww}  M. Pollicott, H. Weiss \& S. A. Wolpert, Topological dynamics of the Weil-Petersson geodesic flow,
{\it Adv. Math.}, {\bf 223} (2010), 1225--1235.

\bibitem{robert} D. Robert, Sur les traces d'op\'erateurs (de Grothendieck \`a Lidskii),  {\it Gaz. Math.}, 
{\bf 141} (2014), 76-91.



\bibitem{ruelleinventiones} D. Ruelle, Zeta-functions 
for expanding maps
  and Anosov flows,  
{\it Invent. Math.}, {\bf 34} (1976), 231--242

\bibitem{ruellebook} D.~Ruelle, {\it Thermodynamic Formalism},
  Reading, Mass., Addison-Wesley, 1978.

\bibitem{ruellerepellers}
D. Ruelle,
Repellers for real analytic maps
{\it Ergod. Th. Dynam. Sys.},
{\bf 2} (1982), 99--107.

\bibitem{ruellefredholm}
D. Ruelle, An extension of the theory of Fredholm determinants,
{\it Publ. Math. (IHES)},
{\bf 72} (1990), 175--193. 

\bibitem{sarnakcmp} P. Sarnak,
Determinants of Laplacians,
{\it Comm. Math. Phys.},
{\bf 110} (1987), 113--120.

\bibitem{shallit}
J. Shallit,
Real numbers with bounded partial quotients: a survey,
{\it Enseign. Math.}, {\bf 38} (1992), 151--187.

\bibitem{smalebullams}
 S. Smale, Differentiable dynamical systems, 
{\it Bull. Amer. Math. Soc.}, {\bf 73} (1967), 747--817. 

\bibitem{smale}
S. Smale, Mathematical problems for the next century, {\it Math. Intelligencer}, {\bf 20} (1998), 7--15.


\bibitem{tuckersmale}
W. Tucker,
A rigorous ODE solver and Smale's 14th problem,
{\it Found. Comp. Math.}, {\bf 2} (2002), 53--117.


\bibitem{ulam}
S. Ulam, {\it Problems in Modern Mathematics},
Interscience, New York, 1960.


\bibitem{Zin00}
M. Zinsmeister,
{\it Thermodynamic formalism and holomorphic dynamical systems}, Translated from the 1996 French original by C. Greg Anderson. SMF/AMS Texts and Monographs, 2. American Mathematical Society, Providence, RI; Soci\'et\'e Math\'ematique de France, Paris, 2000.

\end{thebibliography}
\end{document}